\documentclass[a4paper]{amsproc}
\usepackage[T1]{fontenc}
\usepackage[latin1]{inputenc}
\usepackage[lite]{amsrefs}
\usepackage{amssymb}
\usepackage{amsfonts}
\usepackage{amsxtra}
\usepackage{enumerate}

\DeclareMathOperator{\K}{K}

\DeclareMathOperator{\supp}{supp}

\DeclareMathOperator{\PREC}{\mathrm{Pt}}
\DeclareMathOperator{\VONN}{\mathrm{vN}}

\newcommand*{\Cpl}[1]{{#1}^{\mathrm{c}}}

\DeclareMathOperator{\Hom}{Hom}

\DeclareMathOperator{\coker}{coker}

\newcommand*{\C}{{\mathbb{C}}}
\newcommand*{\R}{{\mathbb{R}}}

\newcommand*{\N}{{\mathbb{N}}}
\newcommand*{\clN}{{\overline{\mathbb{N}}}}

\newcommand*{\Mat}{{\mathbb{M}}}
\newcommand*{\Comp}{{\mathbb{K}}}

\newcommand*{\ID}{{\mathrm{id}}}

\newcommand*{\CONT}{\mathcal{C}}
\newcommand*{\CCINF}{\mathcal{D}}
\newcommand*{\CINF}{\mathcal{E}}

\newcommand*{\specrad}{\rho}

\newcommand*{\brd}{-\hspace{0pt}}
\newcommand*{\nbd}{\nobreakdash-\hspace{0pt}}
\newcommand*{\Cstar}{$C^\ast$\nobreakdash-\hspace{0pt}}

\newcommand*{\abs}[1]{\lvert#1\rvert}
\newcommand*{\norm}[1]{\lVert#1\rVert}

\newcommand*{\hot}{\mathbin{\hat{\otimes}}}
\newcommand*{\prot}{\mathbin{\hat{\otimes}_\pi}}

\newcommand*{\blank}{{\llcorner\!\!\lrcorner}}
\newcommand*{\defeq}{\mathrel{:=}}

\newcommand*{\into}{\rightarrowtail}
\newcommand*{\prto}{\twoheadrightarrow}                                         
\theoremstyle{plain}
\newtheorem{theorem}{Theorem}[section]
\newtheorem{proposition}[theorem]{Proposition}
\newtheorem{lemma}[theorem]{Lemma}
\newtheorem{corollary}[theorem]{Corollary}
\newtheorem{deflemma}[theorem]{Definition and Lemma}

\theoremstyle{definition}
\newtheorem{definition}[theorem]{Definition}
\theoremstyle{remark}
\newtheorem{remark}[theorem]{Remark}

\begin{document}

\title{Bornological versus topological analysis in metrizable spaces}
\author{Ralf Meyer}
\address{Mathematisches Institut\\
         Westfälische Wilhelms-Universität Münster\\
         Einsteinstr.\ 62\\
         48149 Münster\\
         Germany
}
\email{rameyer@math.uni-muenster.de}


\thanks{This research was supported by the EU-Network \emph{Quantum
    Spaces and Noncommutative Geometry} (Contract HPRN-CT-2002-00280)
  and the \emph{Deutsche Forschungsgemeinschaft} (SFB 478).}

\begin{abstract}
  Given a metrizable topological vector space, we can also use its von
  Neumann bornology or its bornology of precompact subsets to do
  analysis.  We show that the bornological and topological approaches
  are equivalent for many problems.  For instance, they yield the same
  concepts of convergence for sequences of points or linear operators,
  of continuity of functions, of completeness and completion.  We also
  show that the bornological and topological versions of
  Grothendieck's approximation property are equivalent for Fréchet
  spaces.  These results are important for applications in
  noncommutative geometry.  Finally, we investigate the class of
  ``smooth'' subalgebras appropriate for local cyclic homology and
  apply some of our results in this context.
\end{abstract}
\maketitle

\section{Introduction}
\label{sec:intro}

Bornological vector spaces provide an ideal setting for many problems
in noncommutative geometry and representation theory.  I met them
first when I studied entire cyclic cohomology in my thesis
(\cite{Meyer:Analytic}).  They are also quite useful for many other
purposes.  They give rise to a very nice theory of smooth
representations of locally compact groups (\cite{Meyer:Smoothrep}).
They allow to take into account the analytical extra structure on
sheaves of smooth or holomorphic functions
(\cite{Schneiders:Quasi-Abelian}).  The projective bornological tensor
product still gives good results for spaces like LF-spaces where the
projective topological tensor product does not work.  This is useful
in order to define cyclic type homology theories for convolution
algebras of non-compact Lie groups and étale groupoids because these
algebras are only LF.  The bornological approach is also very
convenient for spaces of test functions and distributions.

The main motivation for this article is \emph{local cyclic homology},
which is due to Michael Puschnigg (\cite{Puschnigg:Homotopy}).  It is
the first cyclic theory that yields reasonable results for
\Cstar{}algebras.  Puschnigg defines his theory on a rather
complicated category whose objects are inductive systems of ``nice''
Fréchet algebras.  A much more natural setup is the category of
bornological algebras.  However, since most of the analysis that
Puschnigg needs is only worked out for topological vector spaces, he
is forced to use more complicated objects.  Nevertheless, inside the
proofs he treats inductive systems of Fréchet spaces as if they were
bornological vector spaces.

Let~$V$ be a topological vector space.  A subset of~$V$ is called
\emph{von Neumann bounded} if it is absorbed by each neighborhood of
zero.  It is called \emph{precompact} if it can be covered by finitely
many sets of the form $x+U$, $x\in V$, for each neighborhood of
zero~$U$.  The von Neumann bounded and the precompact subsets form two
standard bornologies on~$V$, which we call the \emph{von Neumann
  bornology} and the \emph{precompact bornology} on~$V$, respectively.
In order to define the local cyclic homology of, say, a
\Cstar{}algebra~$A$, we must view~$A$ as a bornological algebra with
respect to the precompact bornology.  We cannot take the von Neumann
bornology because various kinds of approximations can only be made
uniform on precompact subsets.  Thus we replace~$A$ by a rather
complicated bornological algebra and we have to do analysis in~$A$
bornologically.

The main theme of this article is that topological and bornological
analysis in a metrizable topological vector space~$V$ give equivalent
answers to many questions.  Since this observation has its own
intrinsic interest, we analyze some matters in greater depth than
needed for cyclic homology.  We treat both the precompact and the von
Neumann bornology, although we only use the precompact one in
applications.  We do not require convexity unless we really need it.

In the last section we indicate how some of our results apply in
connection with local cyclic homology.  Since the definition of that
theory also involves advanced homological algebra, we do not define it
here.  Nevertheless, we can explain why it is invariant under passage
to ``smooth'' subalgebras.  This is the crucial property of the
theory.  Using the examples of ``smooth'' subalgebras that we exhibit
in Section~\ref{sec:isoradial_examples} this invariance result
implies the good homological properties of local cyclic homology for
\Cstar{}algebras.

We only need a subalgebra to be closed under holomorphic functional
calculus in order to get an isomorphism on topological
$\K$\nbd{}theory.  Since this condition merely involves a single
algebra element, it is an intrinsically commutative concept.  We shall
instead use the \emph{spectral radius} for a bounded subset of a
bornological algebra, which controls the convergence of power series
in several non-commuting variables.  A bounded homomorphism with
``locally dense'' range that preserves the spectral radii of bounded
subsets is called \emph{isoradial}.  This is the concept of ``smooth''
subalgebra that is appropriate for local cyclic homology.  We show
that an isoradial homomorphism is an \emph{approximate local homotopy
equivalence} or more briefly, an \emph{apple}, provided a certain
approximation condition is satisfied.  For instance, the algebra
$\CCINF(M)$ of smooth functions with compact support on a smooth
manifold~$M$ is an isoradial subalgebra of $\CONT_0(M)$ and the
embedding is an apple.  Local cyclic homology is defined so that
apples become isomorphisms in the bivariant local cyclic homology.
Thus it produces equally good results for small algebras like
$\CCINF(M)$ and large algebras like $\CONT_0(M)$.

To study local cyclic homology for bornological algebras, we have to
carry over quite a few analytical concepts known for topological
vector spaces to the bornological setting.  We need continuous and
smooth functions from manifolds into bornological vector spaces,
completeness and the completion, approximation of operators on bounded
subsets, the approximation property and ``dense'' subsets.  We show
that these bornological concepts are equivalent to the corresponding
topological ones if~$V$ is metrizable and given the precompact
bornology.  Along the way we introduce some further important
properties like local separability, bornological metrizability and
subcompleteness.  Many results that hold for the precompact and von
Neumann bornologies on metrizable topological vector spaces extend to
arbitrary bornologically metrizable bornological vector spaces.

It is important for our applications that our definitions and
constructions are \emph{local} in the sense that they are compatible
with \emph{direct unions}.  Let~$(V_i)$ be an inductive system of
bornological vector spaces with injective structure maps and let~$V$
be its direct limit.  Then the natural maps $V_i\to V$ are injective,
so that the~$V_i$ are isomorphic to subspaces of~$V$.  A subset of~$V$
is bounded if and only if it is bounded in one of the subspaces~$V_i$.
Hence it is appropriate to call~$V$ a direct union of the inductive
system~$(V_i)$.  Any separated convex bornological vector space can
be written as a direct union of normed spaces in a canonical way.
Thus if a construction is compatible with direct unions, we can reduce
from the case of separated convex bornological vector spaces to the
case of normed spaces.  This simplifies analysis in convex
bornological vector spaces.

We choose our spaces of continuous and smooth functions to be local in
the sense that a continuous (smooth) function into a direct union
$\varinjlim V_i$ is already a continuous (smooth) function into~$V_i$
for some $i\in I$.  There are alternative definitions that are
non-local.  Similarly, the approximation property has both a global
and a local variant, and the local variant suffices for our
applications.  The only non-local construction that we need is the
\emph{completion}.  Its lack of locality means that we have to replace
it by a ``derived functor'' when we define local cyclic homology.
This derived functor agrees with the completion if and only if the
space in question is \emph{subcomplete}, that is, a subspace of a
complete space.  We obtain some sufficient conditions for
subcompleteness.  They imply that the spaces that we must complete to
compute the local cyclic homology of a Fréchet algebra are
subcomplete, so that the problem with completions usually does not
arise in practice.

\subsection{Some notation}
\label{sec:notation}

We call a subset of a bornological vector space \emph{bounded} if it
belongs to the bornology.  This forces us to call the ``bounded''
subsets of a topological vector space ``von Neumann bounded'' because
we usually prefer the precompact bornology.  For a topological vector
space~$V$ let $\PREC(V)$ and $\VONN(V)$ be the bornological vector
spaces obtained by equipping~$V$ with the precompact bornology and the
von Neumann bornology, respectively.

Everything we do works both for real and complex vector spaces.  To
simplify notation we only formulate results for complex vector
spaces.  We refer to~\cite{Hogbe-Nlend:Bornologies} for the elementary
definitions of bornologies, vector space bornologies, convexity and
separatedness.  We call a subset $S\subseteq V$ of a vector space a
\emph{disk} if it is absolutely convex and satisfies $S=\bigcap_{t>1}
tS$.  We let $V_S=\C\cdot S$ be the linear span of~$S$ equipped with
the semi-norm whose unit ball is~$S$.  The condition $\bigcap_{t>1}
tS=S$ insures that~$S$ is the closed unit ball of~$V_S$.  A subset~$S$
of a bornological vector space is called \emph{circled} if
$\lambda\cdot S\subseteq S$ for all $\lambda\in\C$ with
$\abs{\lambda}\le 1$ and $\bigcap_{t>1} tS=S$.  The \emph{circled
  hull} of a bounded subset in a bornological vector space is again
bounded.  Hence any bornology is generated by circled bounded
subsets.

A \emph{null sequence} is a sequence that converges to~$0$.

We write $\Hom(V,W)$ for the space of bounded linear maps between two
bornological vector spaces and for the space of continuous linear maps
between two topological vector spaces.  These are bornological vector
spaces with the bornologies of equibounded and equicontinuous subsets,
respectively.

\section{Functorial properties of the standard bornologies}
\label{sec:categorical}

We discuss some category theoretic properties of the precompact and
von Neumann bornologies.  We compare continuous and bounded
multilinear maps and the topological and bornological completed tensor
products.  We examine the behavior of our bornologies for direct and
inverse limits and their exactness properties.

\subsection{Boundedness versus continuity}
\label{sec:bounded_continuous}

\begin{lemma}  \label{lem:functors_convex_separated}
  Let~$V$ be a topological vector space.  If~$V$ is locally convex
  then $\PREC(V)$ and $\VONN(V)$ are convex bornological vector
  spaces.  The topological space~$V$ is Hausdorff if and only if
  $\PREC(V)$ and $\VONN(V)$ are separated.
\end{lemma}

Let $V,V_1,\dots,V_n,W$ be topological vector spaces.  Let $f\colon
V_1\times \dots\times V_n\to W$ be a multilinear map.  We want to
compare the notions of (joint) continuity and boundedness for~$f$.

\begin{lemma}  \label{lem:functorial}
  If~$f$ is continuous then $\PREC(f)$ and $\VONN(f)$ are bounded.
  Conversely, if the topological vector spaces $V_1,\dots,V_n$ are
  metrizable then the boundedness of $\PREC(f)$ or $\VONN(f)$ implies
  the continuity of~$f$.
\end{lemma}

\begin{proof}
  If~$f$ is continuous then it is evidently bounded for both
  bornologies.  If~$f$ fails to be continuous and $V_1,\dots,V_n$ are
  metrizable, there are null sequences $(v_{j,k})_{k\in\N}$ in~$V_j$
  such that the sequence $f(v_{1,k},\dots,v_{n,k})$ in~$W$ is
  unbounded.  Since the points of a null sequence form a precompact
  set, $f$ is bounded for neither bornology.
\end{proof}

We call~$V$ \emph{bornological} (or \emph{$\PREC$\nbd{}bornological}) if
a subset that absorbs all von Neumann bounded (or precompact) subsets
is already a neighborhood of the origin.  If~$V$ is bornological then
a family of maps $V\to W$ is equibounded for the von Neumann
bornologies if and only if it is equicontinuous.  That is, there is a
bornological isomorphism $\Hom(V,W)\cong\Hom(\VONN(V),\VONN(W))$.
If~$V$ is $\PREC$\nbd{}bornological then an operator $V\to W$ is
continuous if and only if it is bounded for the precompact
bornologies.  We have a bornological isomorphism
$\Hom(V,W)\cong\Hom(\PREC(V),\VONN(W))$.  However, we cannot replace
$\VONN(W)$ by $\PREC(W)$, this fails already for $V=\C$.  Furthermore,
this discussion does not apply to multilinear maps.

The \emph{complete projective (topological) tensor product}~$\prot$
for complete locally convex topological vector spaces is defined by
its universal property (\cite{Grothendieck:Produits_tensoriels}):
continuous linear maps $V\prot W\to X$ correspond to jointly
continuous bilinear maps $V\times W\to X$ for all complete locally
convex topological vector spaces~$X$.  The same universal property
defines the \emph{complete projective (bornological) tensor
  product}~$\hot$ for complete convex bornological vector spaces
(\cite{Hogbe-Nlend:Completions}).  The following result is proved
already in~\cite{Meyer:Analytic}.  We mention it here for the sake of
completeness:

\begin{theorem}  \label{the:tensor_Frechet}
  The functor~$\PREC$ intertwines the complete projective topological
  and bornological tensor products for Fréchet spaces.  That is, there
  is a natural isomorphism $\PREC(V\prot W) \cong
  \PREC(V)\hot\PREC(W)$ for all Fréchet spaces $V,W$.
\end{theorem}

The corresponding result for the incomplete tensor products also
holds.  It follows from Lemma~\ref{lem:inverse_limits} and
Theorem~\ref{the:metrizable_subcomplete} that
$\PREC(V\mathbin{\otimes_\pi} W)$ and $\PREC(V)\otimes\PREC(W)$ are
both bornological subspaces of the completed tensor product.

\subsection{Direct and inverse limits and exactness}
\label{sec:limits_exact}

Category theory defines inverse and direct limits of ``diagrams'' in a
category. Special cases of inverse limits are direct products and
kernels of maps.  Arbitrary inverse limits in additive categories are
built out of these special cases: the inverse limit of an arbitrary
diagram is the kernel of a map between direct products.  Special cases
of direct limits are direct sums and cokernels of maps.  Arbitrary
direct limits are constructed as the cokernel of a map between direct
sums.

\begin{lemma}  \label{lem:inverse_limits}
  The functors $\PREC$ and~$\VONN$ commute with arbitrary inverse
  limits and with direct sums.
\end{lemma}

\begin{proof}
  It suffices to prove that the functors commute with direct products
  and direct sums and preserve kernels of linear maps.  The latter
  means that if $V\subseteq W$ is a subspace with the subspace
  topology then $\PREC(V)$ and $\VONN(V)$ carry the subspace
  bornologies on~$V$ from $\PREC(W)$ and $\VONN(W)$.  This assertion
  is trivial.  The assertions about direct products and direct sums
  are easy.
\end{proof}

An \emph{LF-space} is a topological vector space which can be written
as a direct limit of a countable strict inductive system of Fréchet
spaces.  Well-known assertions about bounded subsets of LF-spaces
(see~\cite{Treves:Kernels}) amount to the statement that
$$
\PREC(\varinjlim V_n) \cong \varinjlim \PREC(V_n),
\qquad
\VONN(\varinjlim V_n) \cong \varinjlim \VONN(V_n),
$$
if $(V_n)_{n\in\N}$ is a strict inductive system of Fréchet spaces
or LF-spaces.

Neither $\PREC$ nor~$\VONN$ commute with direct limits, in general,
because they do not preserve cokernels.  Let $f\colon V\to W$ be a
continuous linear map between two topological vector spaces.  The
quotient bornology on $\VONN(W)/f(V)$ consists of all images of von
Neumann bounded subsets of~$W$.  It is clear that such subsets are von
Neumann bounded in $W/f(V)$.  The converse need not hold, that is, it
may be impossible to lift von Neumann bounded subsets of $W/f(V)$
to~$W$.  The corresponding assertion for the precompact bornology is
sometimes true:

\begin{theorem}  \label{the:precompact_fully_exact}
  Let~$W$ be a complete metrizable topological vector space and let
  $V\subseteq W$ be a closed subspace.  Then the precompact bornology
  on $W/V$ is the quotient bornology on $\PREC(W)/\PREC(V)$.  A
  diagram $K\to E\to Q$ of complete metrizable topological vector
  spaces is an extension of topological vector spaces if and only if
  $\PREC(K)\to\PREC(E)\to\PREC(Q)$ is an extension of bornological
  vector spaces.
\end{theorem}

\begin{proof}
  Since~$W$ is complete, the quotient $W/V$ is also complete.  Hence
  any precompact subset is contained in a compact subset, so that it
  suffices to lift compact subsets.  Metrizability allows to do this,
  see \cite{Treves:Kernels}*{Lemma~45.1}.
  
  A diagram $K\overset{i}\to E\overset{p}\to Q$ in an additive
  category is an \emph{extension} if $K\cong\ker p$ and $Q\cong\coker
  i$.  In the topological vector space setting, this means that~$i$ is
  a topological isomorphism onto the subspace $i(K)\subseteq E$ with
  the subspace topology and that the induced map $E/i(K)\to Q$ is a
  topological isomorphism with the quotient topology on $E/i(K)$.  A
  similar description is available for bornological vector spaces.
  Suppose first that $K\into E\prto Q$ is an extension of topological
  vector spaces.  Since $\PREC$ preserves kernels, we have $\ker
  (\PREC(p))=\PREC(K)$.  Since the quotient space~$Q$ is separated,
  the subspace $i(K)$ is closed.  Hence the first assertion of the
  theorem yields $\coker (\PREC(i))=\PREC(Q)$.  Thus $\PREC(K) \to
  \PREC(E)\to \PREC(Q)$ is an extension.  Conversely, suppose
  $\PREC(K)\to\PREC(E)\to\PREC(Q)$ to be an extension.  Since the
  points of a null sequence form a precompact set, this implies that
  any sequence in~$K$ that is a null sequence in~$E$ is at least
  bounded in~$K$ and that any null sequence in~$Q$ can be lifted to a
  bounded sequence in~$E$.  Moreover, $K\into E \prto Q$ is exact as a
  sequence of vector spaces.  Using metrizability we can deduce from
  these facts that $K\into E\prto Q$ is a topological extension.
\end{proof}

Theorem~\ref{the:precompact_fully_exact} and
Lemma~\ref{lem:functorial} imply that the functor~$\PREC$ restricted
to complete metrizable topological vector spaces is fully faithful
and fully exact.

\section{Convergence, continuity and smoothness of functions}
\label{sec:continuity}

We define Cauchy and convergent sequences and ``continuous'' functions
in bornological vector spaces.  The appropriate concepts of continuity
are uniform and locally uniform continuity.  Thus we only consider
functions that are defined on metric spaces.  We allow incomplete
spaces because we want to treat Cauchy sequences as uniformly
continuous functions.  These concepts are local in the sense explained
in Section~\ref{sec:intro}.  Hence they can be described easily for
separated convex bornological vector spaces.  A sequence in~$V$
converges or is Cauchy if and only if it converges or is Cauchy in the
usual sense in the normed space~$V_T$ for some bounded disk~$T$.  A
function into~$V$ is locally uniformly continuous if and only if it is
locally uniformly continuous as a function into the normed space~$V_T$
for some bounded disk~$T$.

The main result of this section is
Theorem~\ref{the:continuity_metrizable}, which asserts that the
topological and bornological versions of locally uniform continuity
are equivalent for metrizable topological vector spaces, both for the
precompact and the von Neumann bornology.  This contains the
corresponding assertions about convergent sequences and Cauchy
sequences as special cases.

We can also define $k$~times continuously differentiable functions and
smooth functions from smooth manifolds into separated convex
bornological vector spaces by locality.  We describe these function
spaces as subspaces of spaces of continuous functions.  Hence
Theorem~\ref{the:continuity_metrizable} implies analogous statements
about $k$~times continuously differentiable and smooth functions.

We then show that the fine bornological topology associated to the
precompact or the von Neumann bornology on a metrizable topological
vector space~$V$ is equal to the given topology.  We show that these
bornologies are complete if and only if~$V$ is complete.  All these
results are easy consequences of
Theorem~\ref{the:continuity_metrizable}.

\subsection{Bornological convergence, continuity and
  differentiability}
\label{sec:def_continuity}

\begin{definition}  \label{def:Cauchy_sequence}
  Let~$V$ be a bornological vector space, let $(x_n)_{n\in\N}$ be a
  sequence in~$V$ and let $x_\infty\in V$.  We say that~$(x_n)$
  \emph{converges} to~$x_\infty$ if there exist a circled bounded
  subset $S\subseteq V$ and a null sequence of positive scalars
  $\epsilon=(\epsilon_n)$ such that $x_n\in S$ for all
  $n\in\N\cup\{\infty\}$ and $x_n-x_\infty\in\epsilon_n\,S$ for all
  $n\in\N$.  We call~$(x_n)$ a \emph{Cauchy} sequence if there are $S$
  and~$\epsilon$ as above such that $x_n\in S$ and
  $x_n-x_m\in\epsilon_m\,S$ for all $n,m\in\N$ with $n\ge m$.  If we
  want to specify $S$ or $(S,\epsilon)$ we speak of
  \emph{$S$\nbd{}convergent} and \emph{$S$\nbd{}Cauchy sequences} and of
  \emph{$(S,\epsilon)$\brd{}convergent} and \emph{$(S,\epsilon)$\brd{}Cauchy
    sequences}.
\end{definition}

If~$S$ is even a bounded disk then~$S$ carries a metric from the norm
on~$V_S$.  By definition, $S$\nbd{}Cauchy sequences and
$S$\nbd{}convergent sequences are nothing but Cauchy sequences and
convergent sequences in the metric space~$S$.

\begin{definition}  \label{def:locally_precompact}
  A function $f\colon X\to Y$ between two metric spaces is called
  \emph{uniformly continuous} if for all $\epsilon>0$ there is
  $\delta>0$ such that $d_Y(f(x),f(y))<\epsilon$ for all $x,y\in X$
  with $d_X(x,y)<\delta$.  It is called \emph{locally uniformly
  continuous} if its restriction to any ball of finite radius is
  uniformly continuous.

  A function $w\colon X\times X\to\R_+$ is called a \emph{continuity
  estimator} if it is locally uniformly continuous and satisfies
  $w(x,x)=0$ for all $x\in X$.
  
  A metric space $(X,d)$ is called \emph{locally precompact} if all
  bounded subsets are precompact (that is, totally bounded).
\end{definition}

\begin{definition}  \label{def:continuous_function}
  Let~$V$ be a bornological vector space, $(X,d)$ a metric space and
  $f\colon X\to V$ a function.  We call~$f$ \emph{locally uniformly
    continuous} if there are a circled bounded subset $T\subseteq V$
  and a continuity estimator $w\colon X\times X\to\R_+$ such that
  $f(x)-f(y)\in w(x,y)\cdot T$ for all $x,y\in X$.  If we want to
  specify $T$ or $(T,w)$ we call~$f$ \emph{locally $T$\nbd{}uniformly
    continuous} or \emph{locally $(T,w)$\brd{}uniformly continuous}.
  
  We let $\CONT(X,V)$ be the space of locally uniformly continuous
  functions $X\to V$.  Let $\xi\in X$.  A subset
  $S\subseteq\CONT(X,V)$ is called \emph{locally uniformly continuous}
  if there exist $(T,w)$ as above such that all $f\in S$ are locally
  $(T,w)$\nbd{}uniformly continuous and satisfy $f(\xi)\in T$.  We
  call~$S$ \emph{locally uniformly bounded} if there is $T\subseteq V$
  as above that absorbs $f(Y)$ for each bounded subset $Y\subseteq X$.
\end{definition}

The locally uniformly continuous subsets of $\CONT(X,V)$ and the
locally uniformly bounded subsets both form vector bornologies on
$\CONT(X,V)$.  We call them the bornologies of locally uniform
continuity and locally uniform boundedness, respectively.  We shall
see that the first combines well with precompact bornologies, whereas
the latter combines well with von Neumann bornologies.  Therefore, we
prefer the bornology of locally uniform continuity.

If the metric space~$X$ is bounded then we may drop the qualifier
``locally'' and speak of \emph{uniformly continuous functions} because
locally uniformly continuous functions between precompact metric
spaces are automatically uniformly continuous.  For the same reason,
we may drop the qualifier ``locally uniformly'' and speak of
\emph{continuous functions} if~$X$ is locally compact.  Let~$X$ be a
second countable locally compact space.  A metric on~$X$ is called
\emph{proper} if all bounded subsets are compact.  We equip~$X$ with
any proper metric that defines its topology.  The space $\CONT(X,V)$
of continuous functions $X\to V$ does not depend on the choice of the
metric.

\begin{remark}  \label{rem:sequences}
  Let $\clN=\N\cup\{\infty\}$ be the one-point-compactification
  of~$\N$.  Then a sequence~$(x_n)$ converges towards~$x_\infty$ if
  and only if the function $\clN\ni n\mapsto x_n$ is continuous.
  Equip $\N\subseteq\clN$ with the induced metric.  Then~$\N$ is
  locally precompact, its completion is~$\clN$.  A sequence $(x_n)$ is
  a Cauchy sequence if and only if the function $\N\ni n\mapsto x_n$
  is uniformly continuous.  Thus $\CONT(\N,V)$ and $\CONT(\clN,V)$ are
  the spaces of Cauchy sequences and of convergent sequences in~$V$,
  respectively.  The bornologies of locally uniform continuity consist
  of the $(S,\epsilon)$\brd{}Cauchy and $(S,\epsilon)$\brd{}convergent sets
  of sequences, respectively.
\end{remark}

\begin{lemma}  \label{lem:functions_local}
  Let~$(V_i)$ be an inductive system of bornological vector spaces
  with injective structure maps and let $V\defeq \varinjlim V_i$ be
  its direct union.  Equip all function spaces with the bornologies of
  locally uniform continuity or boundedness.   The spaces
  $\CONT(X,V_i)$ form an inductive system with injective structure
  maps, and $\CONT(X,V)=\varinjlim \CONT(X,V_i)$.  That is, the
  functor $\CONT(X,\blank)$ is local in the sense that it commutes
  with direct unions.
\end{lemma}

\begin{proof}
  Trivial.
\end{proof}

As a result, if~$V$ is a convex bornological vector space then
$\CONT(X,V)$ is the direct union of the spaces $\CONT(X,V_T)$ for the
bounded disks $T\subseteq V$.  The space $\CONT(X,V_T)$ consists
exactly of the locally uniformly continuous functions between the
metric spaces $X$ and~$V_T$.

\begin{remark}  \label{rem:continuity_not_local}
  There is also a notion of locality with respect to the variable~$X$.
  Let $(U_\alpha)$ be an open covering of the space~$X$ such that each
  bounded subset of~$X$ is covered already by finitely
  many~$U_\alpha$.  We may expect that a function~$f$ for which
  $f|_{U_\alpha}$ is locally uniformly continuous for all~$\alpha$ is
  already locally uniformly continuous.  However, this fails with the
  definition above, so that our notion of continuity is not local in
  the variable~$X$.  For instance, if~$\N$ is given the discrete
  metric $d(n,m)\defeq \abs{n-m}$ then a function $f\colon X\to V$ is
  locally uniformly continuous if and only if $f(X)\subseteq \C\cdot
  T$ for some circled bounded subset $T\subseteq V$.  If locally
  uniform continuity were local in~$X$ then any function $\N\to V$
  would have to be locally uniformly continuous.  However, this is
  incompatible with locality in the variable~$V$.  It is easy to
  modify the notion of locally uniform continuity so as to get a
  notion that is local in the variable~$X$ but not in~$V$.  Fix
  $\xi\in X$ and let $B_n(\xi)\subseteq X$ be the set of all $x\in X$
  with $d(x,\xi)\le n$.  Define
  $$
  \tilde\CONT(X,V) \defeq \varprojlim \CONT(B_n(\xi),V).
  $$
  The bornologies of locally uniform continuity and boundedness on
  $\CONT(B_n(\xi),V)$ yield corresponding bornologies on
  $\tilde\CONT(X,V)$.  It depends on the situation whether
  $\CONT(X,V)$ or $\tilde\CONT(X,V)$ is more suitable.  For instance,
  one should use~$\tilde\CONT$ to define continuous group
  representations.  The space that is called $\CINF(G,V)$
  in~\cite{Meyer:Smoothrep} is constructed in the same fashion.  Hence
  we prefer to denote it by $\tilde\CINF(G,V)$ here.
\end{remark}

Next we consider differentiability.  In order to reconstruct a
function from its derivatives we need integration, and integrals are
only defined under some convexity hypothesis.  Therefore, it is
reasonable to restrict to separated convex bornological vector
spaces.  Let~$M$ be a second countable smooth manifold.  For a
topological vector space~$V$ we let $\CONT^k(M,V)$ and $\CINF(M,V)$
be the usual topological vector spaces of $k$~times continuously
differentiable and $\CONT^\infty$\brd{}functions $M\to V$, equipped with
the topology of uniform convergence of derivatives up to order~$k$
(or~$\infty$) on compact subsets of~$M$.  If~$V$ is a normed space, we
equip $\CONT^k(M,V)$ and $\CINF(M,V)$ with two bornologies called the
bornology of locally uniform continuity and boundedness.  The latter
is just the von Neumann bornology.  The first is finer and controls,
in addition, the modulus of continuity of the $k$th derivative.  We do
not have to consider derivatives of lower order because the modulus of
continuity of the $j$th derivative is controlled by the norm of the
$j+1$st derivative.  Thus the two bornologies on $\CINF(M,V)$
coincide.

Any separated convex bornological vector space is a direct union of
normed spaces.  Since the functors $\CONT^k(M,\blank)$ and
$\CINF(M,\blank)$ preserve injectivity of continuous linear maps, we
can define $\CONT^k(M,V)$ and $\CINF(M,V)$ as the direct union of the
spaces $\CONT^k(M,V_T)$ and $\CINF(M,V_T)$, respectively, where~$T$
runs through the bounded disks in~$V$.  The spaces $\CONT^k(M,V)$ and
$\CINF(M,V)$ are local in the same sense as $\CONT(X,V)$ (see
Lemma~\ref{lem:functions_local}).  Since all derivatives of a smooth
function are controlled by the same bounded disk, there is a
difference between smooth and $\CONT^\infty$\brd{}functions in the
bornological case.  This distinction is quite important
in~\cite{Meyer:Smoothrep} and also for local cyclic homology.  As in
Remark~\ref{rem:continuity_not_local}, our function spaces are not
local in~$M$, but there is a variant $\tilde\CINF(M,V)$ that is local
in~$M$ and not in~$V$.

Let $\underline{\infty}\defeq\N$ with the discrete topology and let
$\underline{k}=\{0,1,\dots,k\}$ for $k\in\N$.  The topological space
$X_k\defeq\underline{k}\times TM$ is second countable and locally
compact for all $k\in\clN$.  Let $f\colon M\to V$ be $k$~times
differentiable (or smooth for $k=\infty$).  Its $j$th derivative is a
homogeneous function $TM\to V$ in a natural way.  We define
$X_kf\colon X_k\to V$ by taking the $j$th derivative on $\{j\}\times
TM$.  It is possible to characterize the functions $X_kM\to V$ that
are of the form~$X_kf$ by certain integral equations.  This
construction identifies $\CONT^k(M,V)$ and $\CINF(M,V)$ with certain
closed subspaces of $\CONT(X_k,V)$, respectively.  This works both for
topological and bornological~$V$.  The isomorphism is topological in
the first case and bornological in the latter with respect to either
the bornology of locally uniform continuity or the bornology of
locally uniform boundedness.  Thus we can reduce the study of
$\CONT^k(M,V)$ and $\CINF(M,V)$ to the study of continuous functions.

\subsection{Function spaces for metrizable topological vector spaces}
\label{sec:continuous_metrizable}

Let~$V$ be a topological vector space.  Let $f\colon X\to V$ be a
function.  If $f\in\CONT(X,\PREC(V))$ then $f\in\CONT(X,\VONN(V))$,
and if $f\in\CONT(X,\VONN(V))$ then $f\in\CONT(X,V)$.  Moreover, a
locally uniformly bounded subset of $\CONT(X,\VONN(V))$ is necessarily
von Neumann bounded in $\CONT(X,V)$.  If~$X$ is locally precompact
then a locally uniformly continuous subset of $\CONT(X,\PREC(V))$ is
precompact in $\CONT(X,V)$ (compare this with the Arzelà-Ascoli
Theorem).  These assertions are straightforward to prove and need no
hypothesis on~$V$.  The converse implications hold for metrizable~$V$:

\begin{theorem}  \label{the:continuity_metrizable}
  Let~$V$ be a metrizable topological vector space and let~$X$ be a
  metric space.  The space $\VONN \CONT(X,V)$ is equal to
  $\CONT(X,\VONN V)$ with the bornology of locally uniform
  boundedness.  If~$X$ is locally precompact then
  $\PREC \CONT(X,V)$ is equal to $\CONT(X,\PREC V)$ with
  the bornology of locally uniform continuity.
\end{theorem}

\begin{proof}
  We prove first that a precompact subset~$S$ of $\CONT(X,V)$ is
  locally uniformly continuous in $\CONT(X,\PREC V)$ provided~$X$ is
  locally precompact.  Let~$(U_n)$ be a decreasing sequence of closed
  circled neighborhoods of the origin defining the topology of~$V$.
  Let $p_n\colon V\to\R_+$ be the gauge functional of~$U_n$.  This is
  the homogeneous continuous function with closed unit ball~$U_n$.
  Using that the set~$S$ is precompact in $\CONT(X,V)$ one shows that
  the function
  $$
  w_n(x,y) \defeq \sup {} \{ p_n(f(x)-f(y))\mid f\in S\}
  $$
  on $X\times X$ is a continuity estimator.  Fix a base point $\xi\in
  X$.  There exist constants $1>\delta_n>0$ such that
  $$
  w(x,y)\defeq
  \max {} \{
    \delta_n w_n(x,y)^{1/2},
    \delta_n w_n(x,y)\cdot d(x,\xi),
    \delta_n w_n(x,y)\cdot d(y,\xi)
  \mid n\in\N\}
  $$
  is still a continuity estimator.  Define
  $$
  \alpha\colon S\times X\times X\to V,
  \qquad (f,x,y)\mapsto
  \begin{cases}
    \frac{f(x)-f(y)}{w(x,y)} &\text{for $x\neq y$;} \\
    0 &\text{for $x=y$.}
  \end{cases}
  $$
  Let $S_\xi\defeq \{f(\xi)\mid f\in S\}$ and let $T\defeq
  \alpha(S\times X\times X)\cup S_\xi$.  Let $T^\circ$ be the circled
  hull of~$T$.  We have $f(\xi)\in T$ for all $f\in S$ and
  $f(x)-f(y)\in w(x,y)\cdot T$ for all $f\in S$, $x,y\in X$.  Thus~$S$
  is locally $(T^\circ,w)$\nbd{}uniformly continuous.  It remains to
  prove that~$T$ is precompact.  Then~$T^\circ$ is precompact as well.
  Fix $n\in\N$.  We must cover~$T$ by finitely many sets of the form
  $v+U_n$.  Since~$S_\xi$ is evidently precompact, it suffices to
  cover $\alpha(S\times X\times X)$.  The definition of~$w_n$ implies
  $f(x)-f(y)\in w_n(x,y)\cdot U_n$ for all $f\in S$, $x,y\in X$, so
  that
  $$
  \alpha(f,x,y) \in
  \frac{w_n(x,y)}{w(x,y)} \cdot U_n.
  $$
  By definition of~$w$, we have $w_n\le w$ if
  $d(x,\xi)\ge\delta_n^{-1}$ or $d(y,\xi)\ge\delta_n^{-1}$ or
  $w(x,y)\le\delta_n^2$ because $w_n\le \delta_n^{-2}w^2$.  Hence
  $\alpha(f,x,y)\in U_n$ unless $d(x,\xi),d(y,\xi)\le\delta_n^{-1}$
  and $w(x,y)\ge\delta_n^2$.  Let us restrict attention to the
  subset~$X'$ of triples $(f,x,y)$ satisfying these conditions.  This
  is a bounded subset of $S\times X\times X$ on which~$\alpha$ is
  uniformly continuous.  Since~$X$ is locally precompact, $X'$ and
  hence $\alpha(X')\subseteq V$ is precompact.  Thus~$T$ is
  precompact.  Together with the remarks above the theorem this
  finishes the proof that $\PREC \CONT(X,V)=\CONT(X,\PREC V)$.

  Even without the hypothesis that~$X$ be locally precompact, the same
  argument shows that the set~$T$ above is von Neumann bounded.  Hence
  $\CONT(X,\VONN V)=\CONT(X,V)$ as vector spaces for arbitrary~$X$.
  It remains to prove that a von Neumann bounded subset~$S$ of
  $\CONT(X,V)$ is locally uniformly bounded in $\CONT(X,\VONN V)$.
  By hypothesis,
  $$
  T_n \defeq \{f(x) \mid f\in S,\ x\in X,\ d(x,\xi)\le n\} \subseteq V
  $$
  is von Neumann bounded for each $n\in\N$.  The metrizability of
  $\VONN(V)$, which we prove in Section~\ref{sec:def_metrizable},
  yields a single von Neumann bounded subset $T\subseteq V$ that
  absorbs the sets~$T_n$.  Thus~$S$ is locally uniformly bounded in
  $\CONT(X,\VONN V)$.
\end{proof}

Since we have characterized convergent sequences, Cauchy sequences,
continuously differentiable functions and smooth functions in terms of
locally uniform continuity, we get the following corollaries:

\begin{corollary}  \label{cor:convergence_metrizable}
  Let~$V$ be a metrizable topological vector space, let~$(x_n)$ be
  a sequence in~$V$ and let $x_\infty\in V$.  The following
  assertions are equivalent:
  \begin{enumerate}[(i)]
  \item $(x_n)$ converges towards~$x_\infty$ in the topology of~$V$;

  \item $(x_n)$ converges towards~$x_\infty$ in $\PREC(V)$;

  \item $(x_n)$ converges towards~$x_\infty$ in $\VONN(V)$.
  \end{enumerate}
  An analogous statement holds for Cauchy sequences.
\end{corollary}

\begin{corollary}  \label{cor:smooth_metrizable}
  Let~$V$ be a separated locally convex metrizable topological
  vector space.  Let~$M$ be a second countable smooth manifold.  The
  spaces $\VONN \CONT^k(M,V)$ and $\PREC \CONT^k(M,V)$ are equal to
  $\CONT^k(M,\VONN V)$ with the bornology of locally uniform
  boundedness and $\CONT^k(M,\PREC V)$ with the bornology of locally
  uniform continuity, respectively.  Analogous statements hold for
  smooth functions.
\end{corollary}

The analogous assertions for the variants $\tilde\CONT(X,V)$,
$\tilde\CINF(M,V)$, etc., follow from the results above and
Lemma~\ref{lem:inverse_limits}.

\subsection{The fine bornological topology}
\label{sec:fine_born_topology}

\begin{definition}[\cite{Hogbe-Nlend:Bornologies}]  \label{def:closed}
  A subset~$S$ of a bornological vector space is called \emph{closed}
  if any limit of a convergent sequence in~$S$ lies in~$S$.  The
  closed subsets satisfy the axioms for a topology, which we call the
  \emph{fine bornological topology}.
\end{definition}

Thus we get a closure operation and a notion of dense subset in a
bornological vector space.  Bounded linear maps are continuous for
this topology.  A quotient space $W/V$ is separated if and only if
$V\subseteq W$ is closed (\cite{Hogbe-Nlend:Bornologies}).  However,
the fine bornological topology may be quite pathological: the addition
need not be jointly continuous.

\begin{proposition}  \label{pro:born_topology_metrizable}
  Let~$V$ be a metrizable topological vector space.  Then the fine
  bornological topologies on $\PREC(V)$ and $\VONN(V)$ are equal to
  the given topology.
\end{proposition}

\begin{proof}
  A subset of~$V$ is closed if and only if it is sequentially closed.
  Hence the assertion follows from
  Corollary~\ref{cor:convergence_metrizable}.
\end{proof}

\begin{remark}  \label{rem:product_counterexample}
  The fine bornological topology on $\VONN(V)$ is finer than the given
  topology in general.  There may even be bornologically closed linear
  subspaces that are not topologically closed.  Consider, for
  instance, the product $V\defeq \prod_{i\in I}\C$, where~$I$ is an
  uncountable set.  We equip~$V$ with the product topology and
  bornology.  We think of elements of~$V$ as functions $I\to\C$.  Let
  $V_c\subseteq V$ be the set of all functions with countable support.
  This linear subspace is sequentially closed and hence bornologically
  closed.  However, $V_c$ is dense in~$V$.  The quotient space $V/V_c$
  is a complete convex bornological vector space on which there exist
  no bounded linear functionals.  Any bounded linear functional on
  $V/V_c$ is a continuous linear functional on~$V$ that vanishes
  on~$V_c$ and hence everywhere.
\end{remark}

\subsection{Completeness}
\label{sec:complete}

Recall that we identified the spaces of convergent sequences and
Cauchy sequences with $\CONT(\clN,V)$ and $\CONT(\N,V)$, respectively.
Equip both sequence spaces with the bornology of uniform continuity.

\begin{deflemma}  \label{deflem:complete}
  Let~$V$ be a separated bornological vector space.  Then the
  following conditions are equivalent:
  \begin{enumerate}[(i)]
  \item the map $\CONT(\clN,V)\to\CONT(\N,V)$ is a bornological
    isomorphism;

  \item for any circled bounded subset $S\subseteq V$ and any
    sequence of positive scalars~$\epsilon$, there exist a circled
    bounded subset $T\subseteq V$ and a sequence of positive
    scalars~$\delta$ such that any $(S,\epsilon)$\brd{}Cauchy sequence is
    $(T,\delta)$\brd{}convergent;
    
  \item for any circled bounded subset~$S$ there is a circled bounded
    subset~$T$ such that any $S$\nbd{}Cauchy sequence is
    $T$\nbd{}convergent;

  \item any Cauchy sequence in~$V$ converges and for any circled bounded
    subset $S\subseteq V$ the set of limit points of $S$\nbd{}Cauchy
    sequences is again bounded.

  \end{enumerate}
  \emph{We call~$V$ \emph{complete} if it satisfies these equivalent
    conditions.}
\end{deflemma}

\begin{proof}
  Condition~(ii) just makes explicit the meaning of~(i), so that
  (i)$\iff$(ii).  We show (ii)$\Longrightarrow$(iii).  Fix $S$
  and~$\epsilon$ and find $T$ and~$\delta$ as in~(ii).  Let~$(x_n)$ be
  $S$\nbd{}Cauchy.  Then a subsequence of~$(x_n)$ is
  $(S,\epsilon)$\brd{}Cauchy and hence $(T,\delta)$\nbd{}convergent.
  Therefore, $(x_n)$ itself is
  $(S+T,\epsilon+\delta)$\brd{}convergent.  Thus~(iii) holds.  The
  implication (iii)$\Longrightarrow$(iv) is trivial.  We show
  (iv)$\Longrightarrow$(ii).  This finishes the proof.  Given $S$
  and~$\epsilon$, let~$T$ be the set of all limit points of
  $S$\nbd{}Cauchy sequences.  This set is again circled and bounded.
  Let~$(x_n)$ be $(S,\epsilon)$\brd{}Cauchy.  Let~$x_\infty$ be its
  limit, which exists by~(iv).  For any $m\in\N$ we have
  $\epsilon_m^{-1}(x_{m+n}-x_m)\in S$.  This sequence is in fact
  $(S,\epsilon_{m+n}/\epsilon_m)$\nbd{}Cauchy.  Hence its limit lies
  in~$T$.  This means that $x_\infty-x_m\in\epsilon_m T$.  Thus
  $(x_m)$ is $(T,\epsilon)$\brd{}convergent.
\end{proof}

It is clear that completeness is local, that is, hereditary for direct
unions.  It is also hereditary for arbitrary inverse limits because
closed subspaces and direct products of complete spaces are again
complete.

A disk~$T$ in a bornological vector space~$V$ is called
\emph{complete} if~$V_T$ is a Banach space.  Equivalently, $T$ with
the metric from~$V_T$ is a complete metric space.  If~$T$ is complete
then the limit of any $T$\nbd{}Cauchy sequence is contained in~$T$
again.  If~$V$ is complete then the set of all limit points of
$T$\nbd{}Cauchy sequences is a complete bounded disk.  Therefore, a
convex bornological vector space is complete if and only if any
bounded subset is contained in a complete bounded disk.  This is how
Henri Hogbe-Nlend defines completeness
in~\cite{Hogbe-Nlend:Completions}.

\begin{proposition}  \label{pro:extend_to_completion}
  Let~$V$ be a complete bornological vector space, let $(X,d)$ be a
  metric space and let $(\bar{X},\bar{d})$ be its completion.  Let
  $f\colon X\to V$ be a locally uniformly continuous function.
  Then~$f$ has a unique extension to a locally uniformly continuous
  function $\bar{X}\to V$.  This gives a bornological isomorphism
  $\CONT(X,V)\cong\CONT(\bar{X},V)$ for the bornologies of locally
  uniform continuity and boundedness.
\end{proposition}

\begin{proof}
  Any $x\in\bar{X}$ is the limit of a Cauchy sequence~$(x_n)$ in~$X$.
  By uniform continuity $f(x_n)$ is a Cauchy sequence in~$V$.  It has
  a unique limit because~$V$ is complete.  We define $\bar{f}(x)\defeq
  \lim f(x_n)$.  This does not depend on the choice of the
  sequence~$(x_n)$ because any two sequences converging to~$x$ are
  subsequences of a single convergent sequence.  We have to check
  that~$\bar{f}$ is locally uniformly continuous.  Let~$f$ be locally
  $(S,w)$\brd{}uniformly continuous.  We can extend~$w$ to a
  continuity estimator~$\bar{w}$ on~$\bar{X}$.  The set~$S'$ of all
  limit points of $S$\brd{}Cauchy sequences is again bounded.  So is
  the set~$S''$ of all limit points of $S'$\brd{}Cauchy sequences.  As
  in the proof of Lemma~\ref{deflem:complete} one shows first that
  $f(x)-\bar{f}(y)\in \bar{w}(x,y)\cdot S'$ if $x\in X$, $y\in\bar{X}$
  and then $\bar{f}(x)-\bar{f}(y)\in \bar{w}(x,y)\cdot S''$ for all
  $x,y\in\bar{X}$.  Thus~$\bar{f}$ is locally
  $(S'',\bar{w})$\brd{}uniformly continuous.  This shows that
  $\CONT(X,V)\cong\CONT(\bar{X},V)$.  It is clear that this
  isomorphism is compatible with both standard bornologies.
\end{proof}

The following result shows that bornological completeness is weaker
than topological completeness.

\begin{proposition}  \label{pro:topological_complete}
  Let~$V$ be a Hausdorff topological vector space equipped with the
  von Neumann or the precompact bornology.  Suppose that any
  bornological Cauchy sequence in~$V$ is topologically convergent.
  Then~$V$ is bornologically complete.
\end{proposition}

\begin{proof}
  Let $S\subseteq V$ be circled and bounded.  Then the
  closure~$\bar{S}$ of~$S$ is bounded as well.  Since any
  bornologically convergent sequence is topologically convergent, its
  limit point lies in~$\bar{S}$.  Hence~$V$ satisfies condition~(iv)
  of Definition~\ref{deflem:complete}.
\end{proof}

\begin{theorem}  \label{the:metrizable_complete}
  Let~$V$ be a metrizable topological vector space.  Then the
  following are equivalent:
  \begin{enumerate}[(i)]
  \item $V$ is complete as a topological vector space;

  \item $\PREC(V)$ is bornologically complete;

  \item $\VONN(V)$ is bornologically complete.

  \end{enumerate}
\end{theorem}

\begin{proof}
  We may assume~$V$ to be Hausdorff.  The space~$V$ is complete if and
  only if each Cauchy sequence in~$V$ converges.  The same is true for
  $\PREC(V)$ and $\VONN(V)$ by
  Proposition~\ref{pro:topological_complete}.  Hence the assertion
  follows from Corollary~\ref{cor:convergence_metrizable}.
\end{proof}

\section{Some applications of bornological metrizability}
\label{sec:metrizable}

Metrizability is a global property of a bornological vector space that
encodes some properties of the precompact and the von Neumann
bornologies of metrizable topological vector spaces.  It is a very
useful tool in bornological analysis.  Some applications of
metrizability can be found in~\cite{Meyer:Smoothrep}.  We already used
it in the proof of Theorem~\ref{the:continuity_metrizable}.  The local
version of metrizability is a very weak property because any convex
bornological vector space is locally metrizable.  Since we mainly
consider convex bornologies in applications, this concept may not seem
very useful.  Nevertheless, we take the time to prove the following
structure theorem: a bornological vector space is locally metrizable
if and only if it is a direct union of metrizable topological vector
spaces with the von Neumann bornology.  Local density is the correct
notion of density in connection with approximation problems such as
those in Section~\ref{sec:isoradial_apple}.  We use metrizability and
local separability to show that a subset of a metrizable topological
vector space is locally dense with respect to the precompact bornology
if and only if it is topologically dense.  The same holds for the von
Neumann bornology under a mild additional hypothesis.

The completion~$\Cpl{V}$ of a bornological vector space~$V$ is defined
by a universal property.  Let~$V$ be a metrizable topological vector
space with completion~$\bar{V}$.  We identify the completion of
$\PREC(V)$ with $\PREC(\bar{V})$.  The same holds for the von Neumann
bornology under a mild additional hypothesis.  Even for convex~$V$
the natural map $V\to\Cpl{V}$ need not be injective.  This means that
maps defined on bounded subsets of~$V$ need not extend to~$\Cpl{V}$.
Therefore, we must be very careful with completions when we consider
apples in Section~\ref{sec:isoradial_apple}.  Here we avoid such
problems by requiring our algebras to be complete.  However, to define
local cyclic homology we must pass to analytic tensor algebras and
noncommutative differential forms, so that we must complete tensor
products.  A bornological vector space is called subcomplete if the
map $V\to\Cpl{V}$ is a bornological embedding with locally dense
range.  This is the case where completions are harmless.  We show that
locally separable, bornologically metrizable topological vector spaces
are subcomplete.

\subsection{Metrizability and local metrizability}
\label{sec:def_metrizable}

\begin{definition}  \label{def:born_metrizable}
  A bornological vector space is \emph{(bornologically) metrizable} if
  for any sequence $(S_n)_{n\in\N}$ of bounded subsets there is a
  sequence of positive scalars $(\epsilon_n)_{n\in\N}$ such that
  $$
  \sum \epsilon_n\, S_n \defeq
  \bigcup_{N\in\N} \sum_{n=1}^N \epsilon_n\, S_n
  $$
  is bounded as well.  It is called \emph{locally metrizable} if this
  condition holds for the constant sequence $S_n=S$ for any bounded
  subset~$S$.
\end{definition}

\begin{lemma}  \label{lem:complete_born_metrizable}
  Let~$V$ be a metrizable bornological vector space.  Then~$V$ is
  complete if and only if it satisfies the following strengthening of
  the metrizability condition: for any sequences $(S_n)_{n\in\N}$ of
  bounded subsets there is a sequence of positive scalars
  $(\epsilon_n)_{n\in\N}$ such that the infinite series $\sum_{n\in\N}
  \lambda_n x_n$ converge, where $\lambda_n\in\C$ with
  $\abs{\lambda_n}\le\epsilon_n$, $x_n\in S_n$, and these infinite
  sums form a bounded subset of~$V$.  We denote this bounded subset by
  $\sum^\infty \epsilon_n\, S_n$.

  An analogous statement holds for locally metrizable bornological
  vector spaces.
\end{lemma}

\begin{proof}
  Suppose that~$V$ is complete and let~$(S_n)$ be a sequence of
  circled bounded subsets.  We can choose~$(\epsilon_n)$ such that
  $T\defeq \sum \epsilon_n\cdot n\cdot S_n$ is bounded.  This insures
  that the infinite series $\sum_{n\in\N} \lambda_n x_n$ in the
  statement of the lemma are $T$\nbd{}Cauchy.  Completeness yields
  that they are $U$\nbd{}convergent for some bounded subset~$U$.
  Therefore, $\sum^\infty \epsilon_n\, S_n$ is bounded.

  Suppose conversely that~$V$ satisfies the strengthening of the local
  metrizability condition.  Fix a bounded set~$S$.  Then there is a
  sequence of scalars~$(\epsilon_n)$ such that $T\defeq \sum^\infty
  \epsilon_n\, S$ is well-defined and bounded.  We claim that any
  $S$\nbd{}Cauchy sequence~$(y_n)$ converges towards an element
  of~$T$.  This implies that~$V$ is complete.  We can find a
  subsequence $(y_{n(k)})$ such that $y_{n(k)}-y_{n(k-1)}\in
  \epsilon_k\, S$ for all $k\in\N$.  The claim now follows from
  $y_{n(k)} = \sum_{j=0}^k y_{n(j)}-y_{n(j-1)}$.
\end{proof}

\begin{theorem}  \label{the:metrizable}
  Let~$V$ be a metrizable topological vector space.  Then $\PREC(V)$
  and $\VONN(V)$ are bornologically metrizable.
\end{theorem}

\begin{proof}
  Let $(S_n)_{n\in\N}$ be a sequence of precompact or bounded subsets.
  Let~$(U_n)$ be a decreasing sequence of closed neighborhoods of the
  origin that defines the topology of~$V$.  We may assume that
  $U_{n+1}+U_{n+1}\subseteq U_n$ for all $n\in\N$.  Choose
  $\epsilon_n>0$ such that $\epsilon_n\, S_n\subseteq U_n$.  This
  implies $\sum_{n=m+1}^M \epsilon_n\, S_n \subseteq \sum_{n=m+1}^M
  U_n \subseteq U_m$, using repeatedly that $U_n+U_n\subseteq
  U_{n-1}$.  Hence $\sum \epsilon_n\,S_n\subseteq \sum_{n\le m}
  \epsilon_n\,S_n+U_m$.  The set $\sum_{n\le m} \epsilon_n\,S_n$ is
  precompact or bounded if the sets~$S_n$ are.  Therefore, $\sum
  \epsilon_n\,S_n$ is precompact or bounded, respectively.
\end{proof}

Any convex bornological vector space is locally metrizable because
$\sum \epsilon_n\, S$ is contained in the disked hull of~$S$ once
$\sum \abs{\epsilon_n}\le 1$.  Thus local metrizability is a very weak
condition.

\begin{theorem}  \label{the:locally_metrizable}
  A bornological vector space is locally metrizable if and only if it
  is a direct union of metrizable topological vector spaces equipped
  with the von Neumann bornology.  Analogous statements hold for
  separated or complete locally metrizable spaces: they are direct
  unions of spaces of the form $\VONN(W)$ for separated or complete
  metrizable topological vector spaces~$W$.
\end{theorem}

\begin{proof}
  Local metrizability is evidently hereditary for direct unions and
  $\VONN(V)$ is locally metrizable if~$V$ is a metrizable topological
  vector space.  Therefore, direct unions of metrizable topological
  vector spaces are locally metrizable.  For the converse
  implication we begin with some abstract nonsense which requires no
  hypothesis on~$V$ and which is useful in many similar situations.

  Let~$I$ be the set of all injective bounded maps $f\colon
  \VONN(W)\to V$ where~$W$ is a metrizable topological vector space.
  We say $f\le f'$ if $f=f'\circ\VONN(i)$ for some continuous linear
  map $i\colon W\to W'$, which is necessarily injective.  This is a
  partial order on~$I$.  We claim that~$I$ is directed.  That is, for
  any $f,f'\in I$ there exists $g\in I$ with $f\le g$ and $f'\le g$.
  We can obtain~$g$ from the map
  $$
  (f, f') \colon \VONN(W)\oplus \VONN(W') = \VONN(W\oplus W') \to V
  $$
  by dividing out the kernel of $(f,f')$.  Observe that this
  quotient of a metrizable space is again metrizable.  The spaces
  $\VONN(W)$ form an inductive system indexed by~$I$ with injective
  structure maps.  Hence we can form its direct union $\varinjlim
  \VONN(W)$.  The maps $f\colon \VONN(W)\to V$ give rise to an
  injective bounded linear map $\varinjlim \VONN(W)\to V$.  The
  problem is whether this map is a bornological isomorphism.  We have
  to show that each bounded subset $S\subseteq V$ is the image of a
  bounded subset of $\VONN(W)$ for some $f\in I$.  We construct a
  metrizable topological vector space~$W$, a von Neumann bounded
  subset $T\subseteq W$ and a bounded linear map $f\colon \VONN(W)\to
  V$ such that $S=f(T)$.  Dividing out the kernel of~$f$ we obtain an
  element of~$I$.

  We let~$W$ be the vector space of functions $h\colon S\to\C$ with
  finite support and $f(h)\defeq \sum_{x\in S} h(x)x$.  Let
  $T\subseteq W$ be the set of all characteristic functions of
  singletons $\{x\}\subseteq W$.  Then $f(T)=S$.  We have to equip~$W$
  with a metrizable topology.  Since~$V$ is locally metrizable, there
  is a sequence of positive scalars $\epsilon=(\epsilon_n)$ such that
  $\sum \epsilon_n\, S$ is bounded.  We may assume $\lim \epsilon_n=0$
  and that~$\epsilon$ decreases monotonically.  Let
  $\epsilon^{(0)}\defeq\epsilon$.  We define the derived sequences
  $\epsilon^{(n)}$ for $n\ge1$ recursively by $\epsilon^{(n)}_j\defeq
  \max {} \{\epsilon^{(n-1)}_{2j},\epsilon^{(n-1)}_{2j+1}\}$.  Order
  the points $x_1,\dots,x_n$ in $\supp h$ for $h\in W$ so that the
  sequence $h^\ast_j\defeq \abs{h(x_j)}$ is decreasing and let
  $h^\ast_j=0$ for $j>n$.  Let $U^{(n)} \defeq \{h\in W\mid h^\ast\le
  \epsilon^{(n)}\}$.  Since $\epsilon^{(n)}_j>0$ for all $j\in\N$,
  these are absorbing circled subsets of~$W$.  The
  sequences~$\epsilon^{(n)}$ are constructed so that we have
  $U^{(n)}+U^{(n)}\subseteq 2U^{(n-1)}$.  Hence the sets
  $(U^{(n)}/n)_{n\in\N}$ form the neighborhood basis for a metrizable
  vector space topology on~$W$.  The map~$f$ is bounded for this
  topology because even $f(U^{(0)})$ is bounded.  The set~$T$ is
  clearly von Neumann bounded.  Hence we have constructed the required
  map.

  A similar construction yields the finer results for separated and
  complete locally metrizable spaces.
\end{proof}

\subsection{Locally dense subsets and local separability}
\label{sec:locally_dense}

Let~$V$ be a bornological vector space.

\begin{definition}  \label{def:locally_dense}
  A subset $X\subseteq V$ is \emph{locally dense} if for any circled
  bounded subset $S\subseteq V$ there is a circled bounded subset
  $T\subseteq V$ such that any $v\in S$ is the limit of a
  $T$\nbd{}convergent sequence with entries in $X\cap T$.

  A subset $X\subseteq V$ is \emph{sequentially dense} if any $v\in V$
  is the limit of a convergent sequence with entries in~$X$.
\end{definition}

In general, local density is a stronger requirement than sequential
density and the latter is stronger than density with respect to the
fine bornological topology.

\begin{definition}  \label{def:loc_sep}
  We call~$V$ \emph{locally separable} if for any bounded subset
  $S\subseteq V$ there is a countable subset $A\subseteq S$ and a
  circled bounded subset $T\subseteq V$ containing~$S$ such that any
  point of~$S$ is the limit of a $T$\nbd{}convergent sequence with
  entries in~$A$.
\end{definition}

\begin{proposition}  \label{pro:loc_sep}
  Let~$V$ be a metrizable topological vector space.  The precompact
  bornology on~$V$ is always locally separable.  If~$V$ is separable
  then $\VONN(V)$ is locally separable.
\end{proposition}

\begin{proof}
  Equip~$V$ with a metric that defines its topology and restrict it to
  a circled bounded subset $S\subseteq V$.  Thus~$S$ is a bounded
  metric space and the embedding $S\to V$ is uniformly continuous.
  If~$S$ is precompact, then it contains a dense sequence by
  precompactness.  The same holds for bounded~$S$ provided~$V$ is
  separable.  By Theorem~\ref{the:continuity_metrizable} the map
  $S\to V$ is even $T$\nbd{}uniformly continuous for some precompact or
  von Neumann bounded circled subset $T\subseteq V$, depending on
  whether~$S$ is precompact or not.  Thus Cauchy sequences in the
  metric space~$S$ are mapped to $T$\nbd{}Cauchy sequences.
  Let~$\bar{T}$ be the closure of~$T$.  The proof of
  Proposition~\ref{pro:topological_complete} shows that any
  convergent $T$\nbd{}Cauchy sequence is $\bar{T}$\nbd{}convergent.
  Since~$S$ contains a dense sequence, it follows that~$V$ is locally
  separable.
\end{proof}

\begin{proposition}  \label{pro:loc_met_loc_sep}
  A bornological vector space is locally separable and locally
  metrizable if and only if it is a direct union of separable
  metrizable topological vector spaces with the von Neumann bornology.
\end{proposition}

\begin{proof}
  Since local metrizability and local separability are local
  properties, Theorem~\ref{the:metrizable} and
  Proposition~\ref{pro:loc_sep} imply that direct unions of separable
  metrizable topological vector spaces with the von Neumann bornology
  are locally metrizable and locally separable.  Conversely, let~$V$
  be locally metrizable and locally separable.  For any circled
  bounded subset $S\subseteq V$ there is a circled bounded subset
  $T\subseteq V$ and a countable subset $A\subseteq S$ such that any
  point of~$S$ is the limit of a $T$\nbd{}convergent sequence in~$A$.
  Since~$V$ is locally metrizable, there is a metrizable topological
  vector space~$W$ and an injective bounded linear map $f\colon
  \VONN(W)\to V$ such that $f^{-1}(T)$ and $f^{-1}(S)$ are von Neumann
  bounded subsets of~$W$.  The closed linear span~$W'$ of~$A$ in~$W$
  is a separable, metrizable topological vector space.  Since
  $T$\nbd{}convergent sequences are convergent in~$W$, the
  set~$f^{-1}(S)$ is a bounded subset of~$W'$.  The assertion now
  follows from the abstract nonsense part of the proof of
  Theorem~\ref{the:locally_metrizable}.
\end{proof}

Similarly, a bornological vector space is separated, convex and
locally separable if and only if it is a direct union of separable
normed spaces, and complete, convex and locally separable if and only
if it is a direct union of separable Banach spaces.

\begin{theorem}  \label{the:dense_metrizable_loc_sep}
  Let~$V$ be a metrizable, locally separable bornological vector
  space.  Then a subset $X\subseteq V$ is locally dense if and only if
  it is sequentially dense.
\end{theorem}

\begin{proof}
  It is clear that locally dense subsets are sequentially dense.
  Suppose conversely that~$X$ is sequentially dense.  Let $S\subseteq
  V$ be a circled bounded subset.  Since~$V$ is locally separable,
  there are a countable subset $A\subseteq S$ and a circled bounded
  subset~$S'$ such that any $s\in S$ is the limit of an
  $S'$\nbd{}convergent sequence in~$A$.  Since~$X$ is sequentially
  dense, any $v\in A$ is the limit of a sequence $(x_{v,m})_{m\in\N}$
  in~$X$.  This sequence is $T_v$\nbd{}convergent for some circled
  bounded subset~$T_v$.  By bornological metrizability we can find a
  circled bounded subset~$T$ that contains $2\,S'$ and absorbs the
  sets~$T_v$ for all $v\in A$.  Reparametrizing the sequences
  $(x_{v,m})$, we achieve that they are all $(T,1/n)$\brd{}convergent
  towards~$v$.  Write $s\in A$ as the limit of a $T$\nbd{}convergent
  sequence $(v_n)$ in~$A$.  Then the sequence $(x_{v_n,n})$ is a
  sequence in~$X$ that is $T+T$\brd{}convergent towards~$s$.  Thus~$X$ is
  locally dense in~$V$.
\end{proof}

\begin{theorem}  \label{the:locally_dense_metrizable}
  Let~$V$ be a metrizable topological vector space and let $X\subseteq
  V$ be a subset.  Then the following assertions are equivalent:
  \begin{enumerate}[(i)]
  \item $X$ is locally dense in $\PREC(V)$;

  \item $X$ is dense in~$V$ with respect to the given metrizable
    topology;

  \item $X$ is dense in $\PREC(V)$;

  \item $X$ is dense in $\VONN(V)$.

  \end{enumerate}
  If $\VONN(V)$ is locally separable or if~$V$ is a normed space then
  these conditions are also equivalent to~$X$ being locally dense
  in $\VONN(V)$.
\end{theorem}

\begin{proof}
  Proposition~\ref{pro:born_topology_metrizable} yields the
  equivalence of (ii)--(iv).  Local density evidently implies density.
  If~$X$ is dense then it is sequentially dense for $\PREC(V)$ or
  $\VONN(V)$ by Corollary~\ref{cor:convergence_metrizable}.  Hence
  Theorem~\ref{the:dense_metrizable_loc_sep} implies that~$X$ is
  locally dense in $\PREC(V)$ and locally dense in $\VONN(V)$ if the
  latter bornology is locally separable.  Here we also used
  Theorem~\ref{the:metrizable} and Proposition~\ref{pro:loc_sep}.
  If~$V$ is a normed space then density and sequential density in
  $\VONN(V)$ are equivalent for trivial reasons.
\end{proof}

\subsection{Completions and subcomplete spaces}
\label{sec:completions}

Let~$V$ be a bornological vector space.  Its \emph{completion} is a
complete bornological vector space~$\Cpl{V}$ together with a natural
map $i\colon V\to\Cpl{V}$ such that composition with~$i$ induces an
isomorphism $\Hom(\Cpl{V},W)\cong\Hom(V,W)$ for all complete
bornological vector spaces~$W$.  This universal property
determines~$\Cpl{V}$ uniquely up to isomorphism.  Henri Hogbe-Nlend
constructs completions for convex bornological vector spaces
in~\cite{Hogbe-Nlend:Completions}.  An abstract nonsense argument that
uses that completeness is hereditary for products shows that
completions exist for arbitrary~$V$.  We omit this argument because we
are only interested in the special case of subcomplete spaces,
where we construct the completion explicitly.

\begin{definition}  \label{def:subcomplete}
  A bornological vector space~$V$ is called \emph{subcomplete} if the
  map $V\to\Cpl{V}$ is a bornological embedding with locally dense
  range.
\end{definition}

\begin{proposition}  \label{pro:partial_completion}
  Let $i\colon V\to W$ be a bornological embedding with locally dense
  range.  Then $\Cpl{V}\cong\Cpl{W}$.  Suppose that~$W$ is separated
  and that for any circled bounded subset $S\subseteq V$ there is a
  circled bounded subset $T\subseteq W$ such that~$i$ maps
  $S$\nbd{}Cauchy sequences to convergent sequences with limit in~$T$.
  Then $\Cpl{V}\cong W$.
\end{proposition}

\begin{proof}
  We view $V\subseteq W$ and drop~$i$ from our notation.  We claim
  that any bounded map $f\colon V\to X$ into a complete bornological
  vector space~$X$ extends uniquely to a bounded map $\bar{f}\colon
  W\to X$.  By the universal property of the completion this is
  equivalent to $\Cpl{V}\cong\Cpl{W}$.  Local density yields that for
  any bounded subset $S\subseteq W$ there is a circled bounded
  subset $T\subseteq W$ such that any point in~$S$ is the limit of a
  $T$\nbd{}convergent sequence in~$W$ with entries in~$V$.  Let
  $T'\defeq (T+T)\cap V$.  This is a bounded subset of~$V$ because
  $V\subseteq W$ is a bornological embedding.  A $T$\nbd{}convergent
  sequence with entries in~$V$ is $T'$\nbd{}Cauchy.  Thus any $w\in S$
  is the limit of a $T'$\nbd{}Cauchy sequence~$(v_n)$.  The sequence
  $f(v_n)$ is an $f(T')$\brd{}Cauchy sequence in~$X$.  We define
  $\bar{f}(w)\defeq \lim f(v_n)$.  This is well-defined because~$X$ is
  separated.  The map $\bar{f}\colon W\to X$ is bounded and is the
  only bounded map extending~$f$.  Thus $\Cpl{V}\cong\Cpl{W}$.
  
  Let~$(x_n)$ be an $(S,\epsilon)$\brd{}Cauchy sequence for some
  sequence~$\epsilon$.  Then we can find $x_n'\in T\cap V$ with
  $x_n-x_n'\in \epsilon_n\,T$.  Thus~$(x_n)$ converges if and only
  if~$(x_n')$ converges, and both sequences have the same limit.
  Since $V\subseteq W$ is a bornological embedding, the Cauchy
  condition on~$(x_n)$ implies that~$(x_n')$ is a $U$\brd{}Cauchy
  sequence for some circled bounded subset $U\subseteq V$.  We suppose
  that such Cauchy sequences converge in~$W$ and that their limits
  form a bounded subset.  Thus~$W$ is complete and $\Cpl{V}\cong
  W\cong\Cpl{W}$.
\end{proof}

\begin{theorem}  \label{the:completion_metrizable}
  Let~$V$ be a metrizable topological vector space and let~$\bar{V}$
  be its completion as a topological vector space.  Then $\Cpl{(\PREC
    V)}\cong\PREC(\bar{V})$.  Thus~$V$ is subcomplete.  The
  corresponding assertion for the von Neumann bornology holds if~$V$
  is a normed space or if $\VONN(V)$ is locally separable.
\end{theorem}

\begin{proof}
  Theorem~\ref{the:metrizable_complete} asserts that $\PREC(\bar{V})$
  is complete.  Hence the theorem follows from
  Theorem~\ref{the:locally_dense_metrizable} and
  Proposition~\ref{pro:partial_completion}.
\end{proof}

We want to construct the completion using the same recipe as for
metrizable topological vector spaces.  This works at least for
subcomplete bornological vector spaces.  Recall that $\CONT(\clN,V)$
and $\CONT(\N,V)$ are the spaces of convergent and Cauchy sequences,
respectively.  Let $\CONT_0(\N,V)\subseteq\CONT(\clN,V)$ be the
bornological subspace of null sequences.  Let $i\colon
V\to\CONT(\clN,V)$ send $v\in V$ to the corresponding constant
sequence.  Thus $\CONT(\clN,V)\cong \CONT_0(\N,V)\oplus V$.  Let
$$
\bar{V} \defeq \CONT(\N,V)/\CONT_0(\N,V),
$$
equipped with the quotient bornology and the map $i_\ast\colon
V\to\bar{V}$ induced by~$i$.

\begin{proposition}  \label{pro:subcomplete}
  The following assertions are equivalent for a bornological vector
  space~$V$:
  \begin{enumerate}[(i)]
  \item $V$ is subcomplete;

  \item there exists a bornological embedding $V\to W$ in a complete
    bornological vector space~$W$;

  \item the map $\CONT(\clN,V)\to\CONT(\N,V)$ is a bornological
    embedding;

  \item for any circled bounded subset $S\subseteq V$ there is a
    circled bounded subset $T\subseteq V$ such that any $S$\nbd{}Cauchy
    sequence that converges in~$V$ is already $T$\nbd{}convergent;
    
  \item for any circled bounded subset $S\subseteq V$ the set of all
    limit points of convergent $S$\nbd{}Cauchy sequences is bounded.

  \end{enumerate}
  If~$V$ satisfies these equivalent conditions then
  $\Cpl{V}\cong\bar{V}$.
\end{proposition}

\begin{proof}
  The implication (i)$\Longrightarrow$(ii) is trivial.  The functors
  $\CONT(\clN,\blank)$ and $\CONT(\N,\blank)$ clearly preserve
  bornological embeddings.  Thus (ii) implies~(iii).  The same
  arguments as in the proof of Lemma~\ref{deflem:complete} show that
  (iii)--(v) are equivalent.  Suppose~(iii).  We claim that the map
  $V\to\bar{V}$ satisfies the hypotheses of
  Proposition~\ref{pro:partial_completion}, so that
  $\Cpl{V}\cong\bar{V}$ and~$V$ is subcomplete.  Thus the proof of the
  claim will finish the proof of the proposition.
  
  To check that~$\bar{V}$ is separated, we have to show that
  $\CONT_0(\N,V)$ is closed in $\CONT(\N,V)$.  Since $\CONT_0(\N,V)$
  is a bornological subspace, it suffices to show that if~$(x_n)$ is a
  Cauchy sequence in $\CONT_0(\N,V)$ that converges in $\CONT(\N,V)$
  towards~$x_\infty$, then $x_\infty\in\CONT_0(\N,V)$ as well.  Write
  $x_n=(x_{n,m})_{m\in\N}$, then $x_{\infty,m}=\lim_{n\to\infty}
  x_{n,m}$ for all $m\in\N$.  The Cauchy condition for the sequence
  $(x_n)_{n\in\N}$ easily implies $x_\infty\in\CONT_0(\N,V)$.  The map
  $V\to\bar{V}$ is a bornological embedding because
  $V\cong\CONT(\clN,V)/\CONT_0(\N,V)$ and $\CONT(\clN,V)$ embeds in
  $\CONT(\N,V)$.  Let $S\subseteq V$ be a circled bounded subset and
  let~$\epsilon$ be a null sequence of positive scalars.  Let
  $x=(x_n)$ be an $(S,\epsilon)$\brd{}Cauchy sequence in~$V$.  Then
  $i(x_m)$, $m\in\N$, is a sequence in~$\bar{V}$ that converges
  towards~$[x]$.  To prove this, consider the sequences
  $x^{(m)}_n\defeq x_n$ for $n\le m$ and $x^{(m)}_n=x_m$ for $n\ge m$.
  We have $x^{(m)}-i(x_m)\in\CONT_0(\N,V)$, and $x^{(m)}$ is
  $(T,\sqrt{\epsilon})$\brd{}convergent to~$x$, where
  $T\subseteq\CONT(\N,V)$ is the set of
  $(S,\sqrt{\epsilon})$\brd{}Cauchy sequences.  Therefore, $i(V)$ is
  locally dense in~$\bar{V}$ and Cauchy sequences in~$V$ become
  convergent in~$\bar{V}$ in a controlled fashion.  Thus $V\to\bar{V}$
  satisfies the hypotheses of
  Proposition~\ref{pro:partial_completion}.
\end{proof}

Proposition~\ref{pro:topological_complete} and
Lemma~\ref{lem:inverse_limits} imply that $\PREC(V)$ and $\VONN(V)$
satisfy condition~(ii) of Proposition~\ref{pro:subcomplete} and hence
are subcomplete for any topological vector space~$V$.

\begin{theorem}  \label{the:metrizable_subcomplete}
  A bornological vector space is subcomplete if it is bornologically
  metrizable and locally separable.
\end{theorem}

\begin{proof}
  Let~$V$ be metrizable and locally separable.
  Proposition~\ref{pro:loc_met_loc_sep} allows us to write~$V$ as a
  direct union of an inductive system of separable metrizable
  topological vector spaces $(V_i)_{i\in I}$ equipped with the von
  Neumann bornology.  Consider the inductive system $(\Cpl{V_i})_{i\in
    I}$ of completions.  The structure maps of this system need not be
  injective any more.  For $i\le j$ let $K_{ij}\subseteq\Cpl{V_i}$ be
  the kernel of the map $\Cpl{V_i}\to\Cpl{V_j}$.  Let $K_i\defeq
  \bigcup_{j\ge i} K_{ij}$.  We claim that there is $j\in I_{\ge i}$
  such that $K_i=K_{ij}$.  Before we prove this we show that it
  implies the assertion of the theorem.  The quotients
  $\VONN(\Cpl{V_i}/K_i)$ form an inductive system with injective
  structure maps.  Let~$\bar{V}$ be the direct union of this inductive
  system.  Suppose that $S\subseteq V$ is mapped to a bounded subset
  of~$\bar{V}$.  Then~$S$ is von Neumann bounded in $\Cpl{V_i}/K_i$
  and hence in~$\Cpl{V_j}$ for some $j\ge i$.  Therefore, $S$ is
  already von Neumann bounded in~$V_j$.  Hence the map $V\to\bar{V}$
  is a bornological embedding.  Since the spaces $\Cpl{V_i}/K_i$ are
  complete, $\bar{V}$ is a complete bornological vector space.
  Hence~$V$ is subcomplete.
  
  It remains to find~$j$ with $K_{ij}=K_i$.  Since~$V_i$ is
  separable, so is~$\Cpl{V_i}$.  Hence the subspace
  $K_i\subseteq\Cpl{V_i}$ contains a countable dense subset
  $X\subseteq K_i$.  Elements of~$\Cpl{V_i}$ are limits of Cauchy
  sequences in~$V_i$.  Each $x\in X$ is contained in $K_{ij}$ for
  some $j\in I_{\ge i}$.  We write~$x$ as a limit of a Cauchy sequence
  $(x_n)$ in~$V_i$.  Since $x\mapsto 0$ in~$\Cpl{V_j}$, this Cauchy
  sequence is a null sequence in~$\Cpl{V_j}$ and hence in~$V_j$.
  Therefore, it is $T_x$\nbd{}convergent towards~$0$ for some bounded
  subset $T_x\subseteq V$.  Since~$V$ is metrizable, there is a
  bounded subset $T\subseteq V$ that absorbs the countably many
  subsets~$T_x$ for $x\in X$.  This subset is the image of a von
  Neumann bounded subset of~$V_j$ for some $j\in I_{\ge i}$.  By
  construction, $x\mapsto 0$ in~$V_j$ for all $x\in X$.  Hence the
  closure~$K_i$ of~$X$ in~$\Cpl{V_i}$ is also contained in~$K_{ij}$.
  The inclusion $K_{ij}\subseteq K_i$ is trivial.
\end{proof}

\section{Grothendieck's approximation property}
\label{sec:approximation}

Grothendieck's approximation property is essentially a Banach space
concept.  Hence the extension to convex bornological vector spaces is
just as easy as the extension to locally convex topological vector
spaces.  First we define precompact, compact and relatively compact
subsets and compact operators in the bornological framework.  Then we
explain what kind of approximations of operators we consider.  This is
not quite straightforward because bornological convergence in
$\Hom(V,W)$ usually is too restrictive.  Another issue is that the
Hahn-Banach theorem fails for bornological vector spaces.  It may
happen that there are no globally defined linear functionals.
However, for many applications it is enough to have locally defined
maps.  Hence we consider two variants of the approximation property
which use locally and globally defined linear functionals,
respectively.  They are equivalent for regular spaces.  For Fréchet
spaces the bornological approximation properties for the precompact
and von Neumann bornologies are equivalent to the usual approximation
property in the case of a topological vector space.

\subsection{Compact subsets and compact operators}
\label{sec:precompact_subsets}

\begin{definition}  \label{def:precompact_subset}
  Let~$V$ be a bornological vector space.  A subset $S\subseteq V$ is
  called \emph{(pre)compact} if there is a metric~$d$ on~$S$ such that
  $(S,d)$ is (pre)compact and the map $(S,d)\to V$ is uniformly
  continuous.  It is called \emph{relatively compact} if it is
  contained in a compact subset.
\end{definition}

This definition is local.  That is, a subset of a direct union
$\varinjlim V_i$ is precompact if and only if it is precompact
in~$V_i$ for some $i\in I$, and similarly for compact and relatively
compact subsets.  In particular, if~$V$ is a convex bornological
vector space then a subset is precompact, compact or relatively
compact if and only if it is precompact, etc., in the normed
space~$V_T$ for some bounded disk $T\subseteq V$.

It is easy to see that the precompact and relatively compact subsets
form two vector bornologies on~$V$.  We denote the precompact
bornology on a bornological vector space~$V$ by $\PREC(V)$.  We have
the following implications:
$$
\text{compact} \Longrightarrow
\text{relatively compact} \Longrightarrow
\text{precompact} \Longrightarrow
\text{bounded}.
$$
If~$V$ is complete then precompact$\iff$relatively compact by
Proposition~\ref{pro:extend_to_completion}.

It is often useful to replace a given bornology by the associated
precompact bornology.  This mimics the passage from $\VONN(V)$ to
$\PREC(V)$ in the metrizable case:

\begin{theorem}  \label{the:precompact_metrizable}
  Let~$V$ be a metrizable topological vector space and let $S\subseteq
  V$.  Then the following are equivalent:
  \begin{enumerate}[(i)]
  \item $S$ is topologically precompact;
  \item $S$ is bornologically precompact in $\PREC(V)$;
  \item $S$ is bornologically precompact in $\VONN(V)$.
  \end{enumerate}
  Analogous statements hold for compact and relatively compact
  subsets.
\end{theorem}

\begin{proof}
  The implications (ii)$\Longrightarrow$(iii)$\Longrightarrow$(i) are
  obvious.  To prove (i)$\Longrightarrow$(ii), we equip~$V$ with a
  metric that defines its topology and restrict it to~$S$.  Thus~$S$
  becomes a precompact metric space and the map $S\to V$ is uniformly
  continuous.  Hence it is uniformly continuous as a map to $\PREC(V)$
  by Theorem~\ref{the:continuity_metrizable}.
\end{proof}

\begin{corollary}  \label{cor:precompact_Schwartz}
  Let~$V$ be a metrizable topological vector space.  Any bounded
  subset in $\PREC(V)$ is bornologically precompact.  If~$V$ is
  complete then any bounded subset is bornologically relatively
  compact.
\end{corollary}

\begin{definition}  \label{def:compact_operator}
  Let $V$ and~$W$ be separated convex bornological vector spaces.  An
  operator $f\colon V\to W$ is called \emph{compact} if there is a
  Banach space~$B$ and maps $f_1\colon V\to B$, $f_2\colon B\to W$
  such that $f=f_2\circ f_1$ and~$f_2$ maps the unit ball of~$B$ to a
  compact subset of~$W$.
\end{definition}

Let $V$ and~$W$ be locally convex topological vector spaces.  Suppose
that~$V$ is bornological and that~$W$ is metrizable.  Then an operator
$f\colon \VONN(V)\to\VONN(W)$ is compact if and only if there exists a
neighborhood of the origin~$U$ for which $f(U)$ is relatively compact
in~$W$.  An analogous assertion holds for the precompact bornologies
if~$V$ is $\PREC$\nbd{}bornological.  Hence we get the usual notion of
a compact operator in these cases.

It is not hard to show that the sum of two compact operators is again
compact.  The composition of a compact operator and a bounded operator
(in any order) is again compact.  It is clear that finite rank
operators are compact.  We are mainly interested in the case where~$V$
is a Banach space.  Then an operator $F\colon V\to W$ is compact if
and only if it maps the unit ball of~$V$ to a compact subset of~$W$.
Since the image of the unit ball of~$V$ is automatically complete, it
is compact if and only if it is precompact.

\subsection{Approximation of linear operators}
\label{sec:operator_approximation}

Let $V$ and~$W$ be bornological vector spaces.  Recall that
$\Hom(V,W)$ carries the equibounded bornology.  This gives rise to the
following notion of bornological convergence: a sequence~$(f_n)$ in
$\Hom(V,W)$ converges towards~$f_\infty$ if and only if there exists a
null-sequence $(\epsilon_n)$ and for each bounded subset $S\subseteq
V$ there exists a bounded subset $T\subseteq W$ such that
$(f_n-f_\infty)(S)\subseteq \epsilon_n\, T$ for all $n\in\N$.
However, we usually cannot choose $(\epsilon_n)$ uniformly for
all~$S$.

\begin{definition}  \label{def:local_approximation}
  Let $(f_n)_{n\in\clN}$ be an equibounded family of linear operators
  $V\to W$ and let $S\subseteq V$ be a bounded subset.  We say that
  $(f_n)$ \emph{converges uniformly on~$S$} to~$f_\infty$ if there is
  a bounded subset $T\subseteq W$ and a sequence of scalars
  $(\epsilon_n)$ such that $(f_n-f_\infty)(S)\subseteq\epsilon_n\,T$
  for all $n\in\N$.  We abbreviate this as $(f_n-f_\infty)(S)\to0$.
  We say that $(f_n)$ converges uniformly on bounded, compact or
  precompact subsets if it converges uniformly on all bounded, compact or
  precompact~$S$, respectively.
\end{definition}

Given operators $(f_n)_{n\in\clN}$ we define $F\colon V\to W^\clN$ by
$F(v)(n)\defeq f_n(v)$.  The sequence~$(f_n)$ converges uniformly
on~$S$ to~$f_\infty$ if and only if $F(S)$ is a uniformly continuous
subset of $\CONT(\clN,W)$.  Hence the sequence $(f_n)$ converges
uniformly on bounded subsets if and only if~$F$ is a bounded linear
map $V\to\CONT(\clN,W)$.

\begin{theorem}  \label{the:operator_approximation}
  Let~$V$ be a bornological vector space, let~$W$ be a metrizable
  topological vector space and let $S\subseteq V$ be precompact.  Let
  $(f_n)_{n\in\clN}$ be an equibounded set of linear maps
  $V\to\VONN(W)$.  Then the following are equivalent:
  \begin{enumerate}[(i)]
  \item $(f_n)$ converges towards~$f_\infty$ in the topology of
    uniform convergence on~$S$;

  \item $(f_n-f_\infty)(S)\to0$ in $\PREC(W)$;

  \item $(f_n-f_\infty)(S)\to0$ in $\VONN(W)$;

  \item $\lim f_n(v)=f(v)$ for all $v\in S$.

  \end{enumerate}
  Hence we get the same notion of uniform convergence on (pre)compact
  subsets of~$V$ if we use the topology of~$W$ or the bornologies
  $\PREC(W)$ and $\VONN(W)$.
\end{theorem}

\begin{proof}
  It is clear that
  (ii)$\Longrightarrow$(iii)$\Longrightarrow$(i)$\Longrightarrow$(iv).
  We must prove (iv)$\Longrightarrow$(ii).  Let $\CONT(\clN,W)$ be the
  metrizable topological vector space of continuous functions $\clN\to
  W$.  Define~$F$ as above.  Then~(iv) asserts
  $F(S)\subseteq\CONT(\clN,W)$.  By hypothesis, $S$ is precompact in
  an appropriate topology for which the map $S\to V$ is uniformly
  continuous.  Since the family of operators~$(f_n)$ is equibounded,
  the map $F\colon S\to\CONT(\clN,W)$ is a uniformly continuous map
  between metric spaces.  Hence $F(S)\subseteq \CONT(\clN,W)$ is
  precompact.  By Theorem~\ref{the:continuity_metrizable}, $F(S)$ is
  locally uniformly continuous as a subset of $\CONT(\clN,\PREC(W))$.
  This implies~(ii).
\end{proof}

\begin{definition}  \label{def:approx_finite_rank}
  We say that an operator $f\colon V\to W$ can be \emph{approximated
    uniformly on (pre)compact subsets by finite rank operators} if for
  all (pre)compact subsets $S\subseteq V$ there is a sequence of
  finite rank operators $f_n\colon V\to W$, $n\in\N$, such that
  $(f_n-f)(S)\to0$.
\end{definition}

Definition~\ref{def:approx_finite_rank} allows for a different
sequence of finite rank maps for each precompact subset.  Thus we are
implicitly dealing with a net of operators $f_{S,n}$.  We need nets
already for inseparable Banach spaces.

\subsection{The approximation properties}
\label{sec:approx_property}

Recall that a convex bornological vector space is called
\emph{regular} if the bounded linear functionals on it separate its
points.

\begin{lemma}  \label{lem:finite_rank_map_extension}
  Let $V$ be a regular convex bornological vector space and let~$W$
  be a bornological vector space.  Let $T\subseteq V$ be a bounded
  disk and let $S\subseteq V_T$ be a compact disk.  Let $f\colon
  V_T\to W$ be a bounded finite rank map.  Then there is a sequence
  of bounded finite rank maps $f_n\colon V\to W$, $n\in\N$, such that
  $(f_n)$ converges uniformly on~$S$ to~$f$.
\end{lemma}

\begin{proof}
  We identify finite rank maps $V\to W$ with elements of the
  uncompleted bornological tensor product $V'\otimes W$.  Since the
  map $V_T\to V$ is injective and~$V$ is regular, the image of~$V'$
  in~$V_T'$ is weakly dense.  Since the weak topology and the topology
  of uniform convergence on compact disks have the same continuous
  linear functionals, they also have the same closed convex subsets.
  Hence $V'\subseteq V_T'$ is still dense in the topology of uniform
  convergence on~$S$.  Write $f\in V_T'\otimes V$ as a sum of finitely
  many elementary tensors $l\otimes w$.  For each $l\in V_T'$ there
  is a sequence~$(l_n)$ in~$V'$ that converges towards~$l$ in~$V_S'$.
  Viewing the sum of the elementary tensors $l_n\otimes w$ as a finite
  rank map $V\to W$, we obtain the desired approximation.
\end{proof}

\begin{deflemma}  \label{deflem:local_approximation}
  Let~$V$ be a complete convex bornological vector space.  The
  following conditions are equivalent:
  \begin{enumerate}[(i)]
  \item for any Banach space~$E$ any compact linear map $E\to V$ can
    be approximated uniformly by finite rank operators;

  \item for any Banach space~$E$ any bounded linear map $E\to V$ can
    be approximated uniformly on compact subsets by finite rank
    operators;
    
  \item for any compact disk $S\subseteq V$ there is a bounded disk
    $T\subseteq V$ such that $S\subseteq T$ and the inclusion $V_S\to
    V_T$ is the uniform limit of a sequence of finite rank operators in
    $\Hom(V_S,V_T)$;
    
  \item for any compact disk $S\subseteq V$ there is a compact disk
    $T\subseteq V$ such that $S\subseteq T$ and $V_S\to V_T$ is the
    uniform limit of a sequence of finite rank operators in
    $\Hom(V_S,V_T)$.

  \end{enumerate}
  \emph{If~$V$ satisfies these equivalent conditions we say that~$V$
  has the \emph{local (bornological) approximation property}.}
\end{deflemma}

\begin{proof}
  For a Banach space~$E$ an operator $f\colon E\to V$ is compact if
  and only if it maps the unit ball of~$E$ to a compact disk.  Hence
  condition~(i) holds for all Banach spaces~$E$ and all compact maps
  $E\to V$ once it holds for the inclusions $V_S\to V$ for compact
  disks~$S$.  Condition~(iii) makes explicit what~(i) means for the
  maps $V_S\to V$.  Thus (i)$\iff$(iii).  We next prove the
  implication (i)$\Longrightarrow$(ii).  Again it suffices to
  prove~(ii) for maps of the form $V_T\to V$ for a bounded disk
  $T\subseteq V$.  Let $S\subseteq V_T$ be a compact disk.  By~(i) the
  inclusion $V_S\to V$ can be approximated by finite rank operators
  $V_S\to V$.  Lemma~\ref{lem:finite_rank_map_extension} allows us to
  take bounded finite rank operators $V_T\to V$.  This means that~(ii)
  holds.  It is clear that (iv) implies~(iii).  Hence we are done if
  we prove the implication (ii)$\Longrightarrow$(iv).
  
  Let $S\subseteq V$ be a compact disk.  Thus~$S$ is a compact disk in
  $V_{S''}$ for some complete bounded disk $S''\subseteq V$.  By
  Theorem~\ref{the:precompact_metrizable} there is a compact disk
  $S'\subseteq V_S$ such that~$S$ is already compact in~$V_{S'}$.
  Condition~(ii) provides a sequence of finite rank maps $f_n\colon
  V_{S''}\to V$ such that $(f_n-\ID)(S')\to0$.  This convergence
  happens in~$V_{T'}$ for some bounded disk $T'\subseteq V$.
  Since~$S$ is compact in~$V_{S'}$,
  Theorem~\ref{the:operator_approximation} yields that already
  $(f_n-\ID)(S)\to0$ in $\PREC(V_{T'})$.  That is, $(f_n-\ID)(S)\to0$
  in~$V_T$ for some compact disk $T\subseteq T'$.  This is exactly
  what~(iv) means.
\end{proof}

\begin{deflemma}  \label{deflem:global_approximation}
  Let~$V$ be a complete convex bornological vector space.  The
  following conditions are equivalent:
  \begin{enumerate}[(i)]
  \item the identity map of~$V$ can be approximated uniformly on
    compact subsets by finite rank operators;

  \item any operator $V\to V$ can be approximated uniformly on compact
    subsets by finite rank operators;

  \item for any bornological vector space~$W$ any operator $V\to W$
    can be approximated uniformly on compact subsets by finite rank
    operators;

  \item for any bornological vector space~$W$ any operator $W\to V$
    can be approximated uniformly on compact subsets by finite rank
    operators;

  \item $V$ has the local approximation property and is regular.

  \end{enumerate}
  \emph{If~$V$ satisfies these equivalent conditions we say that~$V$
  has the \emph{global (bornological) approximation property}.}
\end{deflemma}

\begin{proof}
  The equivalence of the first four conditions is proved as for
  topological vector spaces
  (see~\cite{Grothendieck:Produits_tensoriels}).  The idea is that
  $(f_n\circ\phi)$ and $(\phi\circ f_n)$ approximate~$\phi$ on a given
  subset once~$(f_n)$ approximates~$\ID_V$ on a sufficiently large
  subset.  It remains to prove that these conditions are equivalent
  to~(v).  Restricting to Banach spaces in~(iv), we see that the
  global approximation property implies the local one.  If $x\in V$,
  $x\neq0$, then there is a sequence of finite rank maps $f_n\colon
  V\to V$ with $f_n(x)\to x$.  Since finite rank operators are
  elements of $V'\otimes V$, there must be $l\in V'$ with $l(x)\neq0$.
  Hence~(i) implies that~$V$ is regular.  Thus (i)--(iv) imply~(v).

  Conversely, suppose~(v).  Fix a compact disk $S\subseteq V$.
  Then~$S$ is compact in~$V_T$ for some bounded disk $T\subseteq V$.
  By the local approximation property we can approximate the
  inclusion $V_T\to V$ uniformly on~$S$ by bounded finite rank maps
  $V_T\to V$.  Since~$V$ is regular,
  Lemma~\ref{lem:finite_rank_map_extension} allows us to use bounded
  finite rank maps $V\to V$.  Thus~(v) implies~(i).
\end{proof}

The local approximation property is evidently hereditary for direct
unions.

\begin{theorem}  \label{the:approximation_topological}
  Let~$V$ be a Fréchet space.  Then the following are equivalent:
  \begin{enumerate}[(i)]
  \item $V$ has Grothendieck's approximation property as a topological
    vector space;
    
  \item $\VONN(V)$ has the global approximation property;

  \item $\VONN(V)$ has the local approximation property;

  \item $\PREC(V)$ has the global approximation property;

  \item $\PREC(V)$ has the local approximation property.
  \end{enumerate}
\end{theorem}

\begin{proof}
  Since $\PREC(V)$ and $\VONN(V)$ are evidently regular, there is no
  difference between the local and global approximation properties.
  Moreover, $\PREC(V)$ and $\VONN(V)$ have the same compact disks by
  Theorem~\ref{the:precompact_metrizable}.  Since condition~(iv) of
  Definition~\ref{deflem:local_approximation} characterizes the local
  approximation property using only compact disks, the local
  approximation properties for $\VONN(V)$ and $\PREC(V)$ are
  equivalent.  Thus (ii)--(v) are equivalent.  The equivalence
  (i)$\iff$(ii) follows from
  Theorem~\ref{the:operator_approximation}.  The equiboundedness
  requirement in Theorem~\ref{the:operator_approximation} can be
  circumvented as in the proof of the implication
  (ii)$\Longrightarrow$(iv) in Lemma~\ref{deflem:local_approximation}.
\end{proof}

\section{Isoradial homomorphisms and local homotopy equivalences}
\label{sec:isoradial_apple}

Throughout this section, we restrict attention to complete convex
bornological algebras.  The basic concept of this section is the
spectral radius of a bounded subset.  We use it to define locally
multiplicative bornological algebras and isoradial homomorphisms.
Being locally multiplicative means being a direct union of Banach
algebras.  A subalgebra~$A$ of a locally multiplicative algebra~$B$ is
called isoradial if it is locally dense and if a bounded subset of~$A$
has the same spectral radius in $A$ and~$B$.  We exhibit several
important examples of isoradial subalgebras.  Then we introduce
approximate local homotopy equivalences, briefly called apples.  Local
cyclic homology is defined so that apples become isomorphisms in
bivariant local cyclic homology.  Our main theorem asserts that an
isoradial homomorphism is an apple provided a certain approximation
condition is satisfied.  This explains the invariance of local cyclic
homology for ``smooth'' subalgebras and is responsible for the good
properties of the theory.

\subsection{The spectral radius}
\label{sec:specrad}

Let~$A$ be a complete convex bornological algebra.

\begin{definition}  \label{def:spectral_radius}
  Let $S\subseteq A$ be a bounded subset.  We define the
  \emph{spectral radius} $\specrad(S)=\specrad(S;A)$ of~$S$ as the
  infimum of the numbers $r\in\R_{>0}$ for which the set
  $(r^{-1}S)^\infty \defeq \bigcup_{n=1}^\infty (r^{-1}S)^n$ is
  bounded.  If no such~$r$ exists, we put $\specrad(S)=\infty$.  We
  call~$A$ \emph{locally multiplicative} if $\specrad(S)<\infty$ for
  all bounded subsets $S\subseteq A$.
\end{definition}

\begin{proposition}  \label{pro:locally_multiplicative}
  A complete convex bornological algebra is locally multiplicative if
  and only if it is a direct union of Banach algebras.
\end{proposition}

\begin{proof}
  It is clear that direct unions of Banach algebras are locally
  multiplicative.  Suppose conversely that~$A$ is locally
  multiplicative.  Let $S\subseteq A$ be bounded.  Then there is
  $r\in\R_{>0}$ with $\specrad(S)<r$.  The complete disked hull~$T$ of
  $(r^{-1}S)^\infty$ is bounded and satisfies $S\subseteq rT$ and
  $T\cdot T\subseteq T$, so that~$A_T$ is a Banach algebra.  The same
  argument as the abstract nonsense part of the proof of
  Theorem~\ref{the:locally_metrizable} now shows that~$A$ is a direct
  union of Banach algebras.
\end{proof}

The usual Banach algebra functional calculus can be extended easily to
locally multiplicative complete bornological algebras.  In fact, this
was one of the historical motivations to study bornological algebras.

The spectral radius is local in the following sense.  If~$A$ is a
direct union of subalgebras $(A_i)_{i\in I}$ then
\begin{equation}  \label{eq:specrad_union}
  \specrad(S;A) = \liminf \specrad(S;A_i)
\end{equation}
for all bounded subsets $S\subseteq A$.

\begin{lemma}  \label{lem:specrad}
  Let $S\subseteq A$ be a bounded subset.  Let~$S'$ be its disked
  hull.  Then $\specrad(S)=\specrad(S')$.  We have
  $\specrad(cS)=\abs{c}\specrad(S)$ for all $c\in\C$ and
  $\specrad(S^n)=\specrad(S)^n$ for all $n\in\N_{\ge1}$.  Let
  $S_A\subseteq A$ and $S_B\subseteq B$ be bounded disks and let
  $S_A\hot S_B\subseteq A\hot B$ be the complete disked hull of the
  set of elementary tensors $x\otimes y$ with $x\in S_A$, $y\in S_B$.
  Then $\specrad(S_A\hot S_B)\le\specrad(S_A)\cdot\specrad(S_B)$.
\end{lemma}

\begin{proof}
  We only prove $\specrad(S^n)\ge \specrad(S)^n$, the remaining
  assertions are obvious.  Write $S^\infty= \bigcup_{j=0}^{n-1}
  S^j\cdot (S^n)^\infty$.  Hence~$S^\infty$ is bounded once
  $(S^n)^\infty$ is bounded.
\end{proof}

In order to work with the spectral radius, we must have enough subsets
with finite spectral radius.  Therefore, we restrict attention to
locally multiplicative algebras in the following.  However, our
methods still apply in somewhat greater generality.  For instance, the
algebra $\CONT(\R)$ of (unbounded) continuous functions on~$\R$ can
still be treated in a similar way.

\begin{definition}  \label{def:isoradial}
  Let $A$ and~$B$ be locally multiplicative complete convex
  bornological algebras and let $f\colon A\to B$ be a bounded
  homomorphism.  We call~$f$ \emph{isoradial} if $f(A)$ is locally
  dense in~$B$ and $\specrad(S;A)=\specrad(f(S);B)$ for all bounded
  subsets $S\subseteq A$.  If $\ker f=0$, we call~$A$ an
  \emph{isoradial subalgebra} of~$B$.
\end{definition}

\begin{lemma}  \label{lem:isoradial_condition}
  A bounded homomorphism $f\colon A\to B$ with locally dense range is
  isoradial if and only if $\specrad(S;A)\le1$ for all bounded
  $S\subseteq A$ with $\specrad(f(S);B)<1$.
\end{lemma}

\begin{proof}
  Use that $\specrad(f(S))\le\specrad(S)$ always holds and that
  $\specrad(c\,S)=c \specrad(S)$.
\end{proof}

\begin{remark}  \label{rem:locally_multiplicatively_convex}
  Locally multiplicatively convex Fréchet algebras need not be locally
  multiplicative.  For instance, $\PREC\bigl(\prod_{n\in\N} \C\bigr)$
  is not locally multiplicative.  Michael Puschnigg calls a Fréchet
  algebra~$A$ ``nice'' if $\PREC(A)$ is locally multiplicative
  (\cite{Puschnigg:Homotopy}).  For locally multiplicative Fréchet
  algebras our definition of an isoradial subalgebra is equivalent to
  Puschnigg's definition of a smooth subalgebra
  in~\cite{Puschnigg:Homotopy}.
\end{remark}

\begin{theorem}  \label{the:isoradial_stable}
  Let $A$, $B$ and~$C$ be locally multiplicative complete convex
  bornological algebras.  Suppose that~$C$ is nuclear.  If $f\colon
  A\to B$ is an isoradial homomorphism then so is the induced
  homomorphism $f_\ast\colon A\hot C\to B\hot C$.
\end{theorem}

\begin{proof}
  It is clear that $f(A\hot C)$ is locally dense in $B\hot C$.  Let
  $S\subseteq A\hot C$ be a bounded subset with
  $\specrad(f_\ast(S))<1$.  We have to prove $\specrad(S)\le1$.
  Choose~$r$ with $\specrad(f_\ast(S))<r<1$.  Then $T\defeq (r^{-1}
  f_\ast(S))^\infty$ is a bounded subset of $B\hot C$.  Hence~$T$ is
  absorbed by a set of the form $T_B\hot T_C$ with complete bounded
  disks $T_B$ and~$T_C$ in $B$ and~$C$.  Similarly, $S$ itself is
  absorbed by $S_A\hot S_C$ with complete bounded disks $S_A$
  and~$S_C$.  We may assume that~$T_C$ absorbs~$S_C$.  Since~$C$ is
  nuclear, it is a direct union of spaces isomorphic to $\ell^1(\N)$.
  Hence we can choose~$T_C$ such that $C_{T_C}$ is isometric to
  $\ell^1(\N)$.  Since all algebras are locally multiplicative, we may
  assume $S_A$, $T_B$ and~$S_C$ to be submultiplicative and we can
  rescale~$T_C$ so that $\specrad(T_C)\le1$.  By construction we have
  $(r^{-1}f_\ast(S))^\infty\subseteq \beta\cdot T_B\hot T_C$ for some
  $\beta>0$.  Hence
  $$
  f_\ast(S^n)=
  f_\ast(S)^n\subseteq
  r^n\beta\cdot T_B\hot T_C\subseteq
  r\cdot T_B\hot T_C
  $$
  for sufficiently large~$n$.  We fix such an~$n$.  Since $S_A$
  and~$S_C$ are submultiplicative, $S^n$ is still absorbed by $S_A\hot
  S_C$ and hence by $S_A\hot T_C$.  Let $T_\alpha\defeq \alpha\cdot
  S_A\cap f^{-1}(T_B)$ for $\alpha>0$.  This is a bounded disk in~$A$
  with gauge norm
  $$
  \norm{x}_{T_\alpha} =
  \max{} \{ \norm{f(x)}_{T_B}, \alpha^{-1}\norm{x}_{S_A}\} \le
  \norm{f(x)}_{T_B} + \alpha^{-1}\norm{x}_{S_A}.
  $$
  Since $f(T_\alpha)\subseteq T_B$ and~$f$ is isoradial, we
  have $\specrad(T_\alpha)\le1$ and hence $\specrad(T_\alpha\hot
  T_C)\le1$ for all $\alpha>0$.  We want to show that $S^n\subseteq
  T_\alpha\hot T_C$ for sufficiently large~$\alpha$.  Since
  $\specrad(S^n)=\specrad(S)^n$ by Lemma~\ref{lem:specrad}, this
  implies $\specrad(S)\le1$ as desired.
  
  Since~$C_{T_C}$ is isometric to $\ell^1(\N)$, we can estimate the
  gauge norm for $T_\alpha\hot T_C$ as follows.  We have an isometry
  $V\hot C_{T_C}\cong \ell^1(\N,V)$ for any Banach space~$V$.  Hence
  \begin{multline*}
    \norm{x}_{T_\alpha\hot T_C} =
    \norm{x}_{\ell^1(\N,T_\alpha)} =
    \sum_{j\in\N} \norm{x_j}_{T_\alpha}
    \\ \le
    \sum_{j\in\N} \norm{f(x_j)}_{T_B} +  \alpha^{-1}\norm{x_j}_{S_A}
    =
    \norm{f_\ast(x)}_{T_B\hot T_C} +
    \alpha^{-1}\norm{x}_{S_A\hot T_C}.
  \end{multline*}
  For $x\in S^n$ we have $\norm{f_\ast(x)}_{T_B\hot T_C}\le r$ and
  $\norm{x}_{S_A\hot T_C}\le \beta$ for some $\beta>0$.  For
  sufficiently large~$\alpha$ we get $\norm{x}_{T_\alpha\hot
    T_C}\le1$, that is, $S^n\subseteq T_\alpha\hot T_C$.
\end{proof}

\begin{lemma}  \label{lem:PREC_specrad}
  Let~$A$ be a locally multiplicative complete convex bornological
  algebra.  If $S\subseteq A$ is bornologically precompact then
  $\specrad(S;A)=\specrad(S;\PREC A)$.  Thus $\PREC(A)$ is locally
  multiplicative and $\PREC(A)\to A$ is isoradial.
\end{lemma}

\begin{proof}
  Suppose that $(r^{-1}\,S)^\infty$ is bounded.  The lemma follows if
  we show that $(R^{-1}\,S)^\infty$ is precompact for all $R>r$.  Let
  $T'\subseteq A$ be a bounded disk such that~$S$ is precompact
  in~$A_{T'}$.  Let~$T$ be a submultiplicative, complete bounded disk
  that absorbs $(r^{-1}\,S)^\infty\cup T'$.  Hence $(R^{-1}\,S)^n$ is
  precompact in~$A_T$ for all $n\in\N$.  Since~$T$ absorbs
  $(r^{-1}\,S)^\infty$ and $R>r$, for any $\epsilon>0$ there is
  $m\in\N$ such that $(R^{-1}\,S)^n\subseteq \epsilon\,T$ for all
  $n\ge m$.  Thus $(R^{-1}\,S)^\infty$ is precompact in~$A_T$.
\end{proof}

\subsection{Examples of isoradial subalgebras}
\label{sec:isoradial_examples}

Let~$A$ be a locally multiplicative complete bornological algebra.
Let~$M$ be a smooth manifold with countably many connected components
and let $M^+=M\cup\{\infty\}$ be the one point compactification
of~$M$ equipped with any metric that defines its topology.  The space
$\CONT_0(M,A)$ is defined as the subspace of $\CONT(M^+,A)$ of
functions vanishing at~$\infty$.  We equip $\CONT_0(M,A)$ with the
bornology of uniform continuity.  Let $\CCINF(M,A)$ be the space of
smooth compactly supported functions $M\to A$.  This is the direct
union of the spaces $\CINF_0(K,A)$ of smooth functions $M\to A$ with
support in~$K$, where~$K$ runs through the compact subsets of~$M$.

\begin{proposition}  \label{pro:specrad_functions}
  Let~$B$ be $\CONT_0(M,A)$ or $\CCINF(M,A)$ and let $S\subseteq B$ be
  a bounded subset.  For $x\in M$ let $S_x\defeq \{f(x)\mid f\in
  S\}$.  The function $x\mapsto \specrad(S_x;A)$ on~$M$ is upper
  semicontinuous and vanishes at~$\infty$, and
  $$
  \specrad(S;B) = \max{} \{\specrad(S_x;A)\mid x\in M\}.
  $$
\end{proposition}

\begin{proof}
  Write~$A$ as a direct union of Banach algebras~$A_T$.  Then
  $\CONT_0(M,A)$ and $\CCINF(M,A)$ are direct unions of the algebras
  $\CONT_0(M,A_T)$ and $\CCINF(M,A_T)$, respectively.
  By~\eqref{eq:specrad_union} we may assume without loss of generality
  that~$A$ be a Banach algebra.  Let $x\in M^+$ and let
  $r_2>r_1>r_0>\specrad(S_x)$.  Then $(r_0^{-1} S_x)^\infty$ is
  bounded in~$A_T$, so that $(r_1^{-1}S_x)^n\subseteq T$ for
  sufficiently large~$n$.  Since~$S$ is a locally uniformly continuous
  set of functions, we have $(r_2^{-1}S_y)^n\subseteq T$ for~$y$ in
  some neighborhood of~$x$.  Therefore, the function $\specrad(S_x)$
  is upper semicontinuous.  It vanishes at~$\infty$ because
  $S_\infty=\{0\}$.  Therefore, it attains its maximum on~$M$.
  
  Let $r>\specrad(S_x)$ for all $x\in M$.  Then there exist
  $n_j\in\N_{\ge1}$ and an open covering $(U_j)$ of~$M^+$ such that
  $f(y)^n\in r^n\,T$ for all $f\in S$, $y\in U_j$, $n\ge n_j$.
  Since~$M^+$ is compact, we can find a finite subcovering.  Hence we
  can find $n\in\N_{\ge1}$ such that $f(y)^n\in r^n\,T$ for all $y\in
  M^+$, $f\in S$.  This easily implies $\specrad(S^n)\le r^n$ for
  $B=\CONT_0(M,A)$ (use Lemma~\ref{lem:PREC_specrad}).  A
  straightforward computation using the derivation property gives the
  same conclusion for $B=\CCINF(M,A)$ as well.  Hence $\specrad(S)\le
  \max{}\{\specrad(S_x)\}$ as desired.  The converse inequality is
  trivial.
\end{proof}

\begin{proposition}  \label{pro:functions_isoradial}
  The subalgebra $\CCINF(M,A)\subset \CONT_0(M,A)$ is isoradial.
\end{proposition}

\begin{proof}
  The computation in Proposition~\ref{pro:specrad_functions} shows
  that the embedding preserves spectral radii.  To prove that its
  range is locally dense, we can reduce to the case where~$A$ is a
  Banach algebra because both $\CCINF(M,A)$ and $\CONT_0(M,A)$ are
  local in~$A$.  It follows from
  Theorem~\ref{the:locally_dense_metrizable} that $\CCINF(M,A)$ is
  locally dense.
\end{proof}

Theorem~\ref{the:continuity_metrizable} and
Corollary~\ref{cor:smooth_metrizable} yield
$$
\PREC(\CONT_0(M,A))\cong\CONT_0(M,\PREC(A)),
\qquad
\PREC(\CCINF(M,A))\cong \CCINF(M,\PREC(A))
$$
if~$A$ is a Fréchet algebra.  We must use the precompact bornology
because we need the bornology of uniform continuity on $\CONT_0(M,A)$.
In the following \emph{all Fréchet algebras are tacitly equipped with
the precompact bornology}.

\begin{proposition}  \label{pro:union_isoradial}
  Let $(A_i)_{i\in I}$ be an inductive system of locally
  multiplicative complete convex bornological algebras with
  injective structure maps.  Let~$A$ be its direct union and let
  $\iota\colon A\to B$ be an injective bounded homomorphism with
  locally dense range.  Suppose that the composition $A_i\to A\to B$
  is a bornological embedding for all $i\in I$.  Then~$\iota$ is
  isoradial.  The hypotheses above are verified if the~$A_i$ are
  \Cstar{}algebras and~$B$ is the inductive limit \Cstar{}algebra.
\end{proposition}

\begin{proof}
  Since $A_i\to B$ is a bornological embedding, it preserves spectral
  radii.  Any bounded subset of~$A$ is already bounded in~$A_i$ for
  some $i\in I$.  Hence $A\to B$ preserves spectral radii.  Since the
  subalgebra~$A$ is also locally dense in~$B$, it is isoradial.  For
  \Cstar{}algebra inductive limits it is clear that $A_i\to B$
  carries the subspace topology and hence the subspace bornology.  The
  local density of~$A$ follows from
  Theorem~\ref{the:locally_dense_metrizable}.
\end{proof}

Next we consider smooth subalgebras for group actions.  Let $\pi\colon
G\times A\to A$ be a representation of a metrizable locally compact
group~$G$ by automorphisms on a locally multiplicative complete
convex bornological algebra~$A$.  We use the function spaces
$\tilde\CONT(G,A)$ and $\tilde\CINF(G,A)$ defined in
Remark~\ref{rem:continuity_not_local} and~\cite{Meyer:Smoothrep}.
Both are bornological algebras for the pointwise product.  The
representation~$\pi$ is called \emph{locally uniformly continuous} or
\emph{smooth} if $\pi_\ast(g)(a)\defeq \pi(g,a)$ defines a bounded
linear map into $\tilde\CONT(G,A)$ or $\tilde\CINF(G,A)$,
respectively.  The map~$\pi_\ast$ is an algebra homomorphism for the
pointwise product on $\tilde\CONT(G,A)$.  The above notion of
continuous representation is the usual one if~$A$ is a Fréchet algebra
by Theorem~\ref{the:continuity_metrizable} and
Lemma~\ref{lem:functorial}.  The \emph{smooth subspace}~$A_\infty$ for
the group action~$\pi$ is defined in~\cite{Meyer:Smoothrep} as the
intersection
$$
A_\infty \defeq
\tilde\CINF(G,A) \cap \pi_\ast(A)
\subseteq \tilde\CONT(G,A).
$$
It is a closed bornological subalgebra of $\tilde\CINF(G,A)$.  It
is shown in~\cite{Meyer:Smoothrep} that this gives the usual smooth
domain if~$A$ is a Fréchet algebra.

\begin{proposition}  \label{pro:group_isoradial}
  The smooth subalgebra $A_\infty\subseteq A$ for a locally uniformly
  continuous group action is isoradial.
\end{proposition}

\begin{proof}
  In the definition of~$A_\infty$ we can replace $\tilde\CINF(G,A)$
  and $\tilde\CONT(G,A)$ by $\CINF(L,A)$ and $\CONT(L,A)$ for any
  compact neighborhood of the identity $L\subseteq G$
  (see~\cite{Meyer:Smoothrep}).  We claim that the subalgebra
  $\CINF(L,A)\subseteq\CONT(L,A)$ is isoradial.
  Proposition~\ref{pro:specrad_functions} shows that $\CONT(L,A)$ is
  locally multiplicative, so that this assertion makes sense.  In
  order to apply Proposition~\ref{pro:functions_isoradial}, we first
  have to reduce to the Lie group case.  Let $U\subseteq G$ be an
  almost connected open subgroup.  Then we may assume $L\subseteq U$
  and hence can replace~$G$ by~$U$.  Let $k\subseteq U$ be a compact
  normal subgroup for which $U/k$ is a Lie group.  The structure
  theory of almost connected groups yields that~$U$ is the projective
  limit of such quotient groups.  The space $\CINF(L,A)$ is defined as
  the direct union of the spaces $\CINF(L/k,A)$ for such subgroups.
  Since $L/k\subseteq U/k$ is a compact subset of a smooth manifold,
  $\CINF(L/k,A)$ has the usual meaning.  Although $\CONT(L,A)$ is not
  equal to the direct union of the spaces $\CONT(L/k,A)$,
  Proposition~\ref{pro:union_isoradial} yields that $\varinjlim
  \CONT(L/k,A)$ is an isoradial subalgebra of $\CONT(L,A)$.
  Proposition~\ref{pro:functions_isoradial} implies that the
  subalgebras $\CINF(L/k,A)\subseteq\CONT(L/k,A)$ are isoradial.
  Hence $\CINF(L,A)$ is isoradial in $\CONT(L,A)$ as asserted.
  
  The homomorphisms $A_\infty\to\CINF(L,A)$ and $A\to\CONT(L,A)$ are
  bornological embeddings.  Hence they preserve the spectral radii of
  subsets.  Thus the embedding $A_\infty\to A$ preserves spectral
  radii.  It remains to prove that~$A_\infty$ is locally dense in~$A$.
  For any $f\in\CCINF(G)$ convolution with~$f$ defines a bounded
  linear map $A\to A_\infty$.  Explicitly, we have $f\ast a = \int_G
  f(g) \pi(g,a)\,dg$.  If $S\subseteq A$ is bounded then
  $\pi_\ast(S)\subseteq \CONT(L,A)$ is uniformly continuous.  Hence
  the operators of convolution by~$f$ converge to the identity
  uniformly on bounded subsets of~$A$ if~$f$ runs through an
  approximate identity in $\CCINF(G)$.  This implies that~$A_\infty$
  is locally dense in~$A$.
\end{proof}

\subsection{Approximate local homotopy equivalences}
\label{sec:local_homotopy}

In this section we do not want to restrict to locally multiplicative
algebras because the more general case is also important and creates
only minor notational inconveniences.

Let $A$ and~$D$ be separated convex bornological algebras.  Let
$S\subseteq D$ be a bounded disk.  Let $S^{(2)}\subseteq D$ be the
disked hull of $S\cup S\cdot S$.  Let $g\colon D_{S^{(2)}}\to A$ be a
bounded linear map.  Its \emph{curvature} is the bounded bilinear map
$$
\omega_g\colon D_S\times D_S\to A,
\qquad
\omega_g(x,y)\defeq f(xy)-f(x)f(y).
$$
To simplify our notation we write
$$
\abs{g}_\omega \defeq \specrad(\omega_g(S,S);A).
$$
We call~$g$ \emph{approximately multiplicative} if
$\abs{g}_\omega<1$.  (We can replace~$1$ by any $\epsilon>0$ because
$\omega_g(tS,tS)=t^2\omega_g(S,S)$.)  We write $M(S;D,A)$ for the set
of approximately multiplicative maps $D_{S^{(2)}}\to A$.  A
\emph{smooth homotopy} between such maps is an element of
$M(S;D,\CINF([0,1],A))$.  An idea of Joachim Cuntz (\cite{Cuntz:kk})
shows that smooth homotopy is an equivalence relation.  We cannot
directly concatenate smooth homotopies because the derivatives may
jump at the glueing point.  The solution is to reparametrize the
smooth homotopy using a smooth bijection $h\colon [0,1]\to[0,1]$ with
vanishing derivatives at $0$ and~$1$.  We let $H(S;D,A)$ be the set of
smooth homotopy classes of approximately multiplicative maps
$D_{S^{(2)}}\to A$.

Notice that the space $H(S;D,A)$ only depends on things happening in
$D_{S^{(2)}}$.  Hence we may replace~$D$ by the quotient of the tensor
algebra on $D_{S^{(2)}}$ by the ideal generated by the relations
$x\otimes y= x\cdot y$ for $x,y\in S$.  Thus we may restrict attention
to algebras~$D$ with such a ``bounded presentation''.

\begin{definition}  \label{def:local_homotopy}
  Let $f\colon A\to B$ be a bounded homomorphism between two
  separated convex bornological algebras.  We call~$f$ an
  \emph{approximate local homotopy equivalence} or briefly an
  \emph{apple} if the induced map $f_\ast\colon H(S;D,A)\to H(S;D,B)$
  is bijective for any bounded disk~$S$ in any separated convex
  bornological algebra~$D$.
\end{definition}

We can make Definition~\ref{def:local_homotopy} more explicit, but the
result is rather complicated and not particularly useful.  Let
$f\colon A\to B$ be an apple.  Let $T\subseteq B$ be a bounded disk.
Then the inclusion $i_T\colon B_{T^{(2)}}\to B$ defines an element of
$H(T;B,B)$ which must be $f_\ast(g_T)$ for some $g_T\in H(T;B,A)$.  We
represent~$g_T$ by an approximately multiplicative map $g_T\colon
B_{T^{(2)}}\to A$.  These maps play the role of a homotopy inverse
of~$f$.  Since $f_\ast(g_T)=i_T$ in $H(T;B,B)$, there is a smooth
homotopy $h_T^B\in M(T;B,\CINF([0,1],B))$ between $f\circ g_T$
and~$i_T$.  Now let $S\subseteq A$ be a bounded disk.  Then $i_S\colon
A_{S^{(2)}}\to A$ defines an element of $H(S;A,A)$.  The elements
$i_S$ and $g_{f(S)}\circ f\circ i_S\in H(S;A,A)$ are mapped to the
same element of $H(S;A,B)$.  Hence there is a smooth homotopy
$h_S^A\in M(S;A,\CINF([0,1],A))$ between $g_{f(S)} fi_S$ and~$i_S$.
Conversely, the existence of maps $g_T$, $h_T^B$ and $h_S^A$ as above
suffices to guarantee that~$f$ is a local homotopy equivalence.  We
prefer Definition~\ref{def:local_homotopy} because it seems more
tractable.

\begin{theorem}  \label{the:isoradial}
  Let $A$ and~$B$ be locally multiplicative complete convex
  bornological algebras and let $f\colon A\to B$ be an isoradial
  bounded homomorphism.  Suppose that one of the following conditions
  is satisfied:
  \begin{enumerate}[(i)]
  \item any bounded subset of~$B$ is bornologically relatively compact
    and~$B$ has the local approximation property;
    
  \item for each bounded disk $S\subseteq B$ there is a sequence
    $(\sigma_n)$ of bounded linear maps $\sigma_n\colon B_S\to A$ such
    that $(f\circ\sigma_n-\ID)(S)\to0$.

  \end{enumerate}
  Then~$f$ is an approximate local homotopy equivalence (apple).
\end{theorem}

\begin{proof}
  First we claim that condition~(i) implies~(ii).  If~(i) holds then
  any bounded disk $S\subseteq B$ is contained in a compact disk.  By
  the local approximation property there is a complete bounded disk
  $T\subseteq B$ containing~$S$ and a sequence~$(\sigma_n)$ of
  bounded finite rank linear maps $B_S\to B_T$ such that~$(\sigma_n)$
  converges in the norm topology on $\Hom(B_S,B_T)$ towards the
  inclusion map $B_S\to B_T$.  Since $f(A)$ is locally dense, we can
  achieve that $f(A)\cap B_T$ is dense in~$B_T$ by enlarging~$T$.  We
  may replace~$\sigma_n$ by a nearby bounded finite rank map into
  $f(A)\cap B_T$ and lift it to a bounded finite rank map into~$A$.
  The resulting maps verify condition~(ii).  Therefore, we may
  assume~(ii) in the following.
  
  Let~$D$ be a separated convex bornological algebra and let
  $S\subseteq D$ be a bounded disk.  Let $h\colon D_{S^{(2)}}\to B$ be
  a bounded linear map with $\abs{h}_\omega<1$.  Thus $h\in M(S;D,B)$.
  We want to prove that~$h$ is smoothly homotopic to $f\circ h'$ for
  an appropriate $h'\in M(S;D,B)$.  Let~$X$ be the disked hull of
  $h(S^{(2)})+h(S)\cdot h(S)$.  Condition~(ii) yields a sequence of
  bounded linear maps $\sigma_n\colon B_X\to A$ such that
  $(f\circ\sigma_n-\ID)(X)\to0$.  This convergence already happens
  in~$B_{S'}$ for some bounded disk $S'\subseteq B$.  We claim that
  there are $0<r<1$ and a bounded disk $T\subseteq B$ such that
  $T\cdot T\subseteq T$, $\omega_h(S,S)\subseteq rT$ and~$T$ absorbs
  $h(S^{(2)})\cup S'$.
  
  Fix~$R$ with $\abs{h}_\omega<R<1$.  The set $T_1\defeq
  \sqrt{R}\,(R^{-1}\omega_h(S,S))^\infty\subseteq B$ is bounded.  By
  construction, $\omega_h(S,S)\subseteq T_1$ and $T_1\cdot
  T_1\subseteq \sqrt{R}\,T_1$.  Since~$B$ is locally multiplicative,
  there is a submultiplicative bounded disk $T_2\subseteq B$ that
  absorbs the bounded subset $h(S^{(2)})\cup T_1\cup S'$.  Since
  $T_1^2\subseteq \sqrt{R}\,T_1$, we have $T_1^n\subseteq R\,T_2$ for
  sufficiently large~$n$.  Then also $(T_1+\epsilon\, T_2)^n\subseteq
  \sqrt{R}\,T_2$ for some $\epsilon>0$.  The set $T_3\defeq
  T_1+\epsilon\, T_2$ contains $\omega_h(S,S)$, absorbs
  $h(S^{(2)})\cup S'$ and satisfies $\specrad(T_3)<1$.  Finally, the
  disked hull~$T$ of $(r^{-1}T_3)^\infty$ has the required properties
  for any~$r$ between $\specrad(T_3)$ and~$1$.  This establishes the
  claim.
  
  By construction, $B_T$ is a normed algebra.  We obtain a sequence of
  bounded linear operators $h^{(n)}\defeq f\circ\sigma_n\circ h\colon
  D_{S^{(2)}}\to B_T$ that converges uniformly towards~$h$ in
  $\Hom(D_{S^{(2)}},B_T)$.  Since $\omega_h(S,S)\subseteq r\,T$, we
  have $\omega_{h^{(n)}}(S,S)\subseteq \sqrt{r}\,T$ for $n\to\infty$.
  Even more, Proposition~\ref{pro:specrad_functions} yields that for
  $n\to\infty$ the linear homotopy
  $$
  h+t(h^{(n)}-h)\colon D_{S^{(2)}}\to \CINF([0,1],B_T)
  $$
  is approximately multiplicative.  Thus $[h]=[h^{(n)}]$ in
  $H(S;D,B)$.  Since~$f$ is isoradial and $\omega_{f\sigma_n h}(S,S)=
  f(\omega_{\sigma_nh}(S,S))$, we get $\sigma_n h\in M(S;D,A)$.  Thus
  $f_\ast[\sigma_nh]=h$ for sufficiently large~$n$, that is,
  $f_\ast\colon H(S;D,A)\to H(S;D,B)$ is surjective.
  
  Now we prove injectivity.  Let $D\supseteq S$ be as above,
  $h_0,h_1\in M(S;D,A)$ and $H\in M(S;D,\CINF([0,1],B))$ such that
  $H_t=f\circ h_t$ for $t=0,1$.  This means that
  $f_\ast[h_0]=f_\ast[h_1]$ in $H(S;D,B)$.  We have to prove that
  $[h_0]=[h_1]$ in $H(S;D,A)$.  Since the functor
  $\CINF([0,1],\blank)$ is local, $H$ is a bounded map to
  $\CINF([0,1],B_X)$ for some bounded disk $X\subseteq B$.
  Condition~(ii) yields bounded linear maps $\sigma_n\colon
  B_{X^{(2)}}\to A$ such that $(f\circ\sigma_n-\ID)(X^{(2)})\to0$.
  Consider the smooth homotopies $h^{(n)}\colon D_{S^{(2)}}\to
  \CINF([0,1],A)$ defined by $h^{(n)}_t\defeq \sigma_n\circ H_t$ and
  let $H^{(n)}\defeq f_\ast\circ h^{(n)}$.  It is not hard to see that
  $(H^{(n)}-H)(S^{(2)})\to0$.  The induced map
  $$
  f_\ast\colon \CINF([0,1],A) \to \CINF([0,1],B)
  $$
  is isoradial by Proposition~\ref{pro:specrad_functions} or by
  Theorem~\ref{the:isoradial_stable}.  Hence the same argument as in
  the proof of surjectivity shows that $\abs{h^{(n)}}_\omega<1$ for
  $n\to\infty$.  Thus $[h^{(n)}_0]=[h^{(n)}_1]$ in $H(S;D,A)$.  Since
  $f\circ h^{(n)}_0=H^{(n)}_0$ converges uniformly to~$H_0$ for
  $n\to\infty$, the linear homotopy
  $tH_0+(1-t)H^{(n)}_0=f_\ast(th_0+(1-t)h^{(n)}_0)$ is approximately
  multiplicative for $n\to\infty$.  Since~$f_\ast$ is isoradial, we
  get $[h_0]=[h^{(n)}_0]$ for $n\to\infty$.  For the same reason,
  $[h_1]=[h^{(n)}_1]$ and thus $[h_0]=[h_1]$ in $H(S;D,A)$.
\end{proof}

Finally, we examine whether the additional approximation hypothesis of
Theorem~\ref{the:isoradial} holds in the examples in
Section~\ref{sec:isoradial_examples}.  We begin with some general
comments.  Let~$B$ be a locally multiplicative Fréchet algebra.
Theorems \ref{the:precompact_metrizable}
and~\ref{the:approximation_topological} imply that condition~(i) of
Theorem~\ref{the:isoradial} holds if and only if~$B$ has
Grothendieck's approximation property.  In particular, this covers the
case of nuclear \Cstar{}algebras.
Theorem~\ref{the:operator_approximation} yields that the convergence
in~(ii) is equivalent to convergence in the topology of uniform
convergence on~$S$.  The equiboundedness hypothesis in
Theorem~\ref{the:operator_approximation} can be circumvented as in the
proof of Lemma~\ref{deflem:local_approximation}.

Consider now the subalgebra $\CCINF(M,A)\subseteq\CONT_0(M,A)$.  In
this case we can use a sequence of smoothing operators on~$M$ that
converge towards the identity to define maps $\sigma_n\colon
\CONT_0(M,A)\to\CCINF(M,A)$.  It is straightforward to verify that
these maps fulfill condition~(ii) of Theorem~\ref{the:isoradial}.  Thus
the embedding $\CCINF(M,A)\to\CONT_0(M,A)$ is an apple.  For
smoothenings of group representations we have already verified the
approximation condition in the proof of
Proposition~\ref{pro:group_isoradial}.

In the situation of completed direct unions,
Theorem~\ref{the:isoradial} may or may not apply.  Condition~(i) holds
if~$B$ is a Fréchet algebra with Grothendieck's approximation
property.  There are some cases where we can verify condition~(ii)
easily.  For a \Cstar{}algebra~$A$ let $\Comp A\defeq
\Comp(\ell^2\N)\otimes A$ be the \Cstar{}algebra stabilization of~$A$.
This is the \Cstar{}algebra direct limit of the system $(\Mat_nA)$.
Compression to $\C^n\subseteq\ell^2(\N)$ defines maps $\Comp
A\to\Mat_n A$.  Theorem~\ref{the:operator_approximation} shows that
they fulfill condition~(ii) of Theorem~\ref{the:isoradial}.  Hence
$\varinjlim \Mat_n A\to \Comp A$ is an apple.  Similarly, if
$(A_i)_{i\in I}$ is a set of \Cstar{}algebras then the embedding of
the purely algebraic direct sum of $(A_i)$ into the \Cstar{}direct sum
satisfies condition~(ii) of Theorem~\ref{the:isoradial} because we
have bounded projections from the \Cstar{}direct sum onto the factors.

\end{document}